\documentclass[a4paper,12pt,reqno]{amsart}
\usepackage{amssymb}
\usepackage[dvips]{graphicx}
\usepackage{MnSymbol}

\newtheorem{cor}{Corollary}[section]

\DeclareMathOperator{\codim}{codim}  

\input epsf.tex
\newdimen\xsize
\newdimen\oldbaselineskip
\newdimen\oldlineskiplimit
\def\restorelineskip{\baselineskip=\oldbaselineskip%
\lineskiplimit=\oldlineskiplimit}
\def\putm[#1][#2]#3{
\hbox{\vbox to 0pt{\parindent=0pt%
\vskip#2\xsize\hbox to0pt{\hskip#1\xsize $#3$\hss}\vss}}}%
\long\def\Line#1{\hbox to \hsize{#1}}
\def\putt[#1][#2]#3{
\vbox to 0pt{\noindent\hskip#1\xsize\lower#2\xsize%
\vtop{\restorelineskip#3}\vss}}

\makeatletter
\def\xbig[#1]#2{{\hbox{$\m@th\left#2\vbox to#1\xsize{}%
\right.\n@space$}}}
\makeatother
\def\xlar[#1]#2{%
\smash{\mathop{ \hbox to #1\xsize{\leftarrowfill}}\limits^{#2}}}
\def\xrar[#1]#2{%
\smash{\mathop{ \hbox to #1\xsize{\rightarrowfill}}\limits^{#2}}}
\def\xline[#1]{\hbox to #1\xsize{\leaders\hrule\hfill}}

\DeclareFontFamily{U}{rsf}{\skewchar\font'177}%
\DeclareFontShape{U}{rsf}{m}{n}{<-6>rsfs5<6-8>rsfs7<8->rsfs10}{}%
\DeclareFontShape{U}{rsf}{b}{n}{<-6>rsfs5<6-8>rsfs7<8->rsfs10}{}%
\DeclareMathAlphabet\RSFS{U}{rsf}{m}{n}
\SetMathAlphabet\RSFS{bold}{U}{rsf}{b}{n}
\def\sf#1{{\mathsf{#1}}}

\def\slsf{\slshape \sffamily }

\def\msmall#1{\mathchoice{\hbox{\small$\displaystyle {#1}$}}{#1}{#1}{#1}}

\def\ss{{\mathbb S}}

\def\cc{{\mathbb C}}
\def\dd{{\mathbb D}}

\def\rr{{\mathbb R}}

\def\nn{{\mathbb N}}

\def\pp{{\mathbb P}}

\def\adyn{\sf{1}}

\def\c{\sf{c}}

\def\const{\sf{const}}

\def\codim{\sf{codim}\,}

\def\deg{\sf{deg}\,}

\def\dim{\sf{dim}\,}

\def\ker{\sf{Ker}\,}
\def\lim{\mathop{\sf{lim}}}

\def\max{\sf{max}}
\def\merto{\dashedrightarrow}
\def\min{\sf{min}}

\def\Reg{\sf{Reg}\,}

\def\trdeg{\sf{tr.deg}\,}

\def\sup{\sf{sup}\,}

\def\eps{\varepsilon}

\def\<{\langle}\let\la=\<
\def\>{\rangle}\let\ra=\>

\def\d{\partial}

\def\ddef{\mathrel{{=}\raise0.3pt\hbox{:}}}
\def\deff{\mathrel{\raise0.3pt\hbox{\rm:}{=}}}

\def\fraction#1/#2{\mathchoice{{\msmall{ #1\over#2}}}%
{{ #1\over #2 }}{{#1/#2}}{{#1/#2}}}
\def\norm#1{\left\Vert{#1}\right\Vert}

\def\le{\leqslant}

\def\emptyset{\varnothing}

\def\longpoints{\leaders\hbox to 0.5em{\hss.\hss}\hfill \hskip0pt}
\def\stateskip{\smallskip}
\def\state#1. {\stateskip\noindent{\bf#1. }} 
\def\statep#1. {\stateskip\noindent{\bf#1 }} 
\def\proof{\state Proof. \2}

\def\Chi{\raise 2pt\hbox{$\chi$}}

\def\ie{\hskip1pt plus1pt{\sl i.e.\/,\ \hskip1pt plus1pt}}

\def\sli{{\sl i)} } 
\def\slii{{\sl i$\!$i)} } 
\def\sliii{{\sl i$\!$i$\!$i)} }

\def\reg{^\sf{reg}}

\def\Chi{\raise 2pt\hbox{$\chi$}}

\let\phI=\phi\let\phi=\varphi\let\varphi=\phI
\let\cal=\mathcal

\def\calc{{\cal C}}

\def\calh{{\cal H}}

\def\calj{{\cal J}}
\def\calk{{\cal K}}

\def\calm{{\cal M}}

\def\calo{{\cal O}}
\def\calp{{\cal P}}

\def\calu{{\cal U}}
\def\calv{{\cal V}}

\def\calx{{\cal X}}
\def\caly{{\cal Y}}

\def\eps{\varepsilon}

\def\d{\partial}

\def\1{{1\mkern-5mu{\rom l}}}

\def\ge{\geqslant}

\def\fraction#1/#2{\mathchoice{{\msmall{ #1\over#2}}}%
{{ #1\over #2 }}{{#1/#2}}{{#1/#2}}}

\def\le{\leqslant}

\def\emptyset{\varnothing}
\newcommand{\2}{\thinspace}

\def\qed{\ \ \hfill\hbox to .1pt{}\hfill\hbox to .1pt{}\hfill $\square$\par}

\def\comment#1\endcomment{}

 \belowdisplayskip=\abovedisplayskip

\def\lineeqqno(#1){\hfill\llap{\vbox to 10pt%
{\vss\begin{align} \eqqno(#1)\end{align}\vss}}\vskip1pt}

\advance\headheight 1.2pt

\def\ShowwLLabel#1{}

\def\thechpt{\Roman{chpt}}

\def\newchapt[#1]#2{\newpage%
\refstepcounter{chpt}\setcounter{subsection}{0}%
\setcounter{thm}{0}\setcounter{defi}{0}%
\setcounter{rema}{0}\setcounter{exrc}{0}%
\renewcommand{\thesubsection}{\thechpt.\arabic{subsection}}%
\section*{\begin{center}\huge \bf Chapter \thechpt\\
#2 \end{center}}\label{#1}%
\ \smallskip%
\markboth{Chapter \thechpt}{#2}%
}
\def\newsect[#1]#2{\refstepcounter{section}\setcounter{equation}{0}%
\renewcommand{\thesubsection}{\arabic{section}.\arabic{subsection}}%
\section*{\arabic{section}.
#2}\vspace{-20pt}\label{#1}\vspace{20pt}%
\markboth{Section \arabic{section}}{#2}}

\def\newlect[#1]#2{\refstepcounter{section}%
\renewcommand{\thesubsection}{\arabic{section}.\arabic{subsection}}%
\section*{Lecture \arabic{section}\\
#2}\label{#1}%
\markboth{Lecture \arabic{section}}{#2}}

\def\newprg[#1]#2{\refstepcounter{subsection}%
\subsection*{{\thesubsection.\ #2}} \label{#1}%
}

\def\newappx[#1]#2{%
\refstepcounter{appx}\setcounter{section}{0}%
\renewcommand{\thesubsection}{A\arabic{appx}.\arabic{subsection}}%
\section*{Appendix \arabic{appx}\\ #2}
\label{#1}%
\markboth{Appendix A\arabic{appx}}{#2}
}

\newtheorem{mainthm}{Main Theorem}
   \def\newthm#1{\begin{mainthm}}
   \def\refthm#1{{\sl Main Theorem }}

\newtheorem{thm}{Theorem}[section]
   \def\newthm#1{\begin{thm}\label{#1}}
   \def\refthm#1{{\sl Theorem \ref{#1}}}

\newtheorem{nnthm}{Theorem}
   \def\newthm#1{\begin{nnthm}\label{#1}}
   \def\refthm#1{{\sl Theorem \ref{#1}}}
\newtheorem{lem}{Lemma}[section]
   \def\newlemma#1{\begin{lem} \label{#1}}

\newtheorem{prop}{Proposition}[section]
   \def\newprop#1{\begin{prop}\label{#1}}

\newtheorem{nnprop}{Proposition}
   \def\newprop#1{\begin{nnprop}\label{#1}}

\newtheorem{corol}{Corollary}[section]
   \def\newcorol#1{\begin{corol} \label{#1}}

\newtheorem{nncorol}{Corollary}
   \def\newcorol#1{\begin{nncorol} \label{#1}}

\newtheorem{defi}{Definition}[section]
   \def\newdefi#1{\begin{defi} \label{#1}\rm }

\newtheorem{nndefi}{Definition}
   \def\newdefi#1{\begin{nndefi} \label{#1}\rm }

\newtheorem{exmp}{Example}[section]
   \def\newexmp#1{\begin{exmp} \label{#1}\rm }

\newtheorem{exrc}{Exercise}
   \def\newexrc#1{\begin{exrc} \label{#1}\rm }

\newtheorem{rema}{Remark}[section]
   \def\newrema#1{\begin{rema} \label{#1}\rm }

\newtheorem{nnrema}{Remark}
   \def\newrema#1{\begin{nnrema} \label{#1}\rm }

\def\eqqno(#1){\label{(#1)}}
\def\eqqref(#1){(\ref{(#1)})}

\def\el2{\sf{L^2}}
\def\l{\it{l}}
\def\vect{{\mathrm{v}}}

\title{Non-existence of a holomorphic embedding of\\ the Sobolev loop space 
into the \\ projective Hilbert space}
\author{Anakkar M., Ivashkovych S*.}
\address{*Universit\'e de Lille-1, UFR de Math\'ematiques, 59655 Villeneuve
d'Ascq, France.}
\email{serge.ivashkovych@univ-lille.fr}
\email{anakkarmohammed@yahoo.fr}

\subjclass{Primary - 32D15, Secondary - 32A20, 46G20, 46T25}
\keywords{Hilbert manifold, meromorphic map, loop space.}
\thanks{* The second author was partially supported by the  
R-CDP-24-004-C2EMPI project.}
\date{\today}

\begin{document}

\begin{abstract}
The goal of this paper is to understand the properties of meromorphic mappings 
with values in two model complex  Hibert manifolds: projective Hilbert space 
$\pp(l^2)$ and Sobolev loop space of the Riemann sphere $L\pp^1$. It occurs 
that these properties are quite different. Based on our study we obtain as a 
corollary that $L\pp^1$ is not biholomorphic to a submanifold of 
$\pp(l^2)$. In other words $L\pp^1$ is {\slsf not} a Hilbert projective variety
despite of the fact that it is K\"ahler and meromorphic functions separate 
points on it. Moreover, we prove that $L\pp^1$ doesn't admit even a non-degenerate
meromorphic map to $\pp (l^2)$.
\end{abstract}

\maketitle

\setcounter{tocdepth}{1}
\tableofcontents

\newsect[INT]{Introduction}

\newprg[INT.non-imb]{On the embedding of a loop space into a projective space}
By a complex Hilbert manifold we understand a Hausdorff topological space $\calx$ 
locally homeomorphic to open subsets in a separable Hilbert space $L$, \ie $\cc^n$ 
or $\l^2$, such that the transition maps 
are holomorphic. The notation $\l^2$ stands for the Hilbert space of sequences of complex
numbers $z=\{z_k\}_{k=1}^{\infty}$ such that $||z||^2:= \sum_k|z_k|^2<\infty $ with the 
standard Hermitian scalar product $(z,w) = \sum_kz_k\bar w_k$ and the standard basis 
$\{ e_1,e_2,...\}$. By $B(z_0,r)$  
we denote the ball of radius $r$ centred at $z_0$ in $L$, by $B(r)$ the ball centred at 
the origin and by $B$ the unit ball. If we want to underline that we are speaking about 
unit ball in $\l^2$ we write $B^\infty$. We say that meromorphic functions separate 
points on the complex  Hilbert manifold $\calx$ if for every pair $p\not=q$ of points 
of $\calx $ there exists a meromorphic function $f$ on $\calx$ such that $f$ is 
holomorphic near $p$ and $q$ and $f(p)\not=f(q)$. It is not difficult to see that 
a compact ($\Rightarrow$ finite dimensional) complex manifold $X$ such that
meromorphic functions separate points on it admits a {\slsf meromorphic injection} 
to the complex projective space $\pp^N$ for some $N$, see Proposition \ref{mer-inj}.
By a meromorphic injection we mean a meromorphic mapping $f:X\merto \pp^N$ such 
that it is a holomorphic injection outside of its indeterminacy set. 

\smallskip  This turns out to not be the case in infinite dimensions. Namely, despite 
the fact that the loop space $L\pp^1$ of the Riemann sphere possesses the property of 
meromorphic separation, see Proposition \ref{mer-sep}, it cannot be meromorphically 
injected into the projective Hilbert space $\pp (\l^2)$. Even more, the following is true.

\begin{nnthm}
\label{non-imb}
There doesn't exist a meromorphic map $g: L\pp^1 \dashedrightarrow \pp(l^2)$ such that the 
codimension
of $\ker dg_q$ is $\ge 2$ at some point $q$ where $g$ is holomorphic.
\end{nnthm}

Mappings $g$ with $\codim\ker dg_z \le 1$ for all $z$ where $g$ is holomorphic, are
mappings to curves and such obviously do exist.
The non-existence of an {\slsf equivariant holomorphic embedding} with respect to the
natural group action of $\Omega \left(PGL(2,\cc)\right)$ and (any) its representation
in $U(l^2)$ was proved in \cite{Z} responding (partially and negatively) to the statement
from \cite{MZ}. The proof in \cite{Z} is based on the comparison of Picard groups and 
concerns a slightly more general settings of Banach manifolds. Our proof is based on the 
comparison of the extension properties of the meromorphic mappings with values in
$L\pp^1$ and $\pp (l^2)$ and doesn't requires $g$ to be equivariant in any sense.

\newprg[INT.def]{Projective Hilbert space}

\smallskip 
Let $\calx$ and $\caly$ be complex Hilbert manifolds. Intuitively a meromorphic mapping 
$\calx \merto \caly$ should be defined by a holomorphic map $f: \calx \setminus I \to 
\caly$, where $I$ is a Hilbert analytic set in $\calx$ of codimension $\geqslant 2$ such
that the graph $\Gamma_f$ of $f$ ``extends'' in some sense to $\calx\times \caly$. And
this ``extended object'' is what we call {\slsf a meromorphic 
graph} (or, a {\slsf a meromorphic map}) and denote it still as $\Gamma_f$ (or, as $f$).
The indeterminacy set of a meromorphic mapping $f$ is the smallest closed subset $I$ of 
$\calx$ such that $f$ is holomorphic on $\calx\setminus I$. It will be denoted as $I_f$. 
We do not claim that $I_f$ is Hilbert analytic in general, but this will be so in both 
model examples of Hilbert manifolds mentioned above. Furthermore we shall investigate 
in what sense the aforementioned extension is realized in the examples of $L\pp^1$ and 
$\pp (\l^2)$. They are quite different in nature but have one common feature, namely the 
notion of meromorphicity for them is  somehow ``predefined'', as it will be 
explained below. 

\smallskip  We start with the case when $\caly $ is the projective Hilbert space 
$\pp (\l^2) \deff (l^2\setminus\{0\})/w\sim \lambda w$, $\lambda \not= 0$, where
$w=(w_0,w_1,...)\in \l^2$. By $z^j_k = w_j/w_k$ we denote the affine coordinates in
the affine chart $\calu_k = \{[w]:w_k\not=0\} = \pp (l^2)\setminus \calh_k$, where 
$\calh_k$ is the hyperplane 
\begin{equation}
\eqqno(hyp-hj)
\calh_k\deff \{[w]\in \pp (\l^2): w_k=0\}. 
\end{equation}

Let a holomorphic mapping $f:\calx\setminus I \to \pp (l^2)$ be given. Assuming that $f(\calx\setminus I)\not\subset \calh_0$ set $A_0\deff f^{-1}(\calh_0)$. 
Observe that $f|_{\calx\setminus A_0}$ is a holomorphic map with values in $\calu_0\equiv 
\l^2$ and therefore writes as $f=(f_1,f_2,...)$. Here $f_k$ are holomorphic functions on 
$\calx\setminus A_0$ and such that $\norm{f(z)}_{\l^2}^2 = \sum_k|f_k(z)|^2$ is locally 
bounded, see subsection \ref{MER-F.hol}. We prove the following statement. 
\begin{nnthm}
\label{pl2-mer}
The closure $\bar A_0$ of the divisor $A_0$ is a divisor in $\calx$ and all $f_k$ extend 
meromorphically to the whole of $\calx$. Moreover, the orders of the poles of these extensions
are locally uniformly bounded along $\bar A_0$. In addition the following is true. 
 
\sli For every point $a\in \calx$ there exists a neighbourhood $\calv\ni a$ and holomorphic 
in  $\calv$  
 
\quad functions $\phi_0,\phi_1,...$ without common factors such that 
 
\smallskip\quad a) $\norm{\phi(z)}^2 = \sum_k|\phi_k(z)|^2$ is locally bounded and $f_k = 
\phi_k/\phi_0$  on $\calv$;
 
\smallskip\quad  b) $\calv\cap I\supset \{z\in \calv : \phi_k(z)=0 \text{ for all } k\}$.

\slii  The indeterminacy set $I_f$ of the meromorphic map $f$ thus obtained  is locally 

\quad defined by the equation $\phi (z) = 0$, \ie $I_f\cap \calv = \{z\in \calv:\phi_0(z)=
\phi_1(z) =...0\}$.
\end{nnthm}

Now it is clear that by a meromorphic mapping in the case of $\pp (\l^2)$ we should understand 
the object obtained by the extension from $\calx\setminus I$ to $\calx$ of a holomorphic map $f$
as in this theorem. It is worth to note that on $\calv\setminus I_f$ our map is
defined as 
\begin{equation}
\eqqno(loc-mer)
f(z) = [\phi_0(z):\phi_1(z):...]
\end{equation}
for an appropriate holomorphic map $\phi = (\phi_0,\phi_1,...)$ with values in $\l^2$.
Therefore it is natural to give the following definition.
\begin{nndefi} 
\label{mer-map1}
A meromorphic mapping from a Hilbert manifold $\calx$ to the complex projective Hilbert 
space $\pp (\l^2)$ is defined by a sequence $f_1,f_2,...$ of meromorphic functions on 
$\calx$ such that locally there exist holomorphic functions $\phi_0,\phi_1,...$  with 
$\norm{\phi}_{\l^2}$ bounded and $f_k=\frac{\phi_k}{\phi_0}$ for all $k\in \nn$. 
\end{nndefi}

\smallskip Notice that such mapping $f$ is locally given by \eqqref(loc-mer) outside of 
the set of common zeroes of $\phi_0,\phi_1,...$. We prove Thullen, Levi and a Hartogs 
type extension theorems for meromorphic mappings with values in $\pp (l^2)$, see Theorems 
\ref{levi-thul} and \ref{H-ext} as well as Corollary \ref{concave}. In this Introduction 
we shall need the following corollary from the Hartogs type extension.

\begin{nncorol}
\label{hart-corol}
Let $\Gamma$ be a $\calc^1$-smooth closed curve in $\cc^2$ and let $f:\cc^2\setminus 
\Gamma \to \pp (\l^2)$ be a holomorphic map. Then $f$ extends meromorphically to the 
whole of $\cc^2$.
\end{nncorol}
This corollary follows from the Hartogs type extension theorem by placing a Hartogs 
figure $H$ near a given point $a\in \Gamma$ in such a way that $H\cap \Gamma = \emptyset$ 
but the associated to $H$ polydisk $\dd^2$ contains $a$.

\newprg[LPS]{Case of loop spaces.}
Our second model example when the meromorphicity is "predefined" is the Sobolev loop 
space of the Riemann sphere $\pp^1$, i.e. $L\pp^1= W^{1,2}(\ss^1,\pp^1)$, where the 
latter denotes the space of maps $u:\ss^1\to \pp^1$ which are of Sobolev class 
$W^{1,2}$. The natural complex structure on $L\pp^1$ was introduced 
by Lempert, see \cite{L3}, and will be recalled in section \ref{LOOP}. For the moment, 
it will be sufficient to say that a holomorphic map $f:\calx \to L\pp^1$ from a complex 
Hilbert manifold $\calx$ to $L\pp^1$ can be seen as a 
mapping $f: \calx \times \ss^1 \to \pp^1$ such that: 

\begin{itemize}
\item for every $s \in \ss^1$ the map $f(\cdot, s): \calx \to \pp^1$ is holomorphic, 
\item for every $z \in \calx$ one has $f(z,\cdot) \in W^{1,2}(\ss^1,\pp^1)$ and the 
correspondence 
\[
\calx \ni z \mapsto f(z,\cdot) \in W^{1,2}(\ss^1, \pp^1) 
\]
is continuous with respect to the standard topology on $\calx$ and the Sobolev 
topology on $W^{1,2}(\ss^1, \pp^1)$. 
\end{itemize} 
Such description will be called the {\slsf representation} of $f$. We shall prove 
the following

\begin{nnthm}
\label{lp1-mer}
Any holomorphic mapping $f:\calx\setminus I\to L\pp^1$, where $\calx$ is a Hilbert manifold
and $I$ a Hilbert analytic set in $\calx$ of codimension $\ge 2$ possesses the following properties:

\sli for every $s\in \ss^1$ the map $f(\cdot ,s)$ is a meromorphic function which meromorphically 

\quad extends  to the whole of $\calx$;
 
\slii the family of indeterminacy sets $\{I_{f(\cdot , s)}:s\in \ss^1\}$ is locally finite 
in $\calx$;
 
\sliii mapping $f$ holomorphically extends to $\calx\setminus \bigcup_{s\in \ss^1}
 I_{f(\cdot ,s)}$, and therefore $I_f = \bigcup_{s\in \ss^1}
 I_{f(\cdot ,s)}$ 
 
\quad is the indeterminacy set of the meromorphic mapping thus obtained.
\end{nnthm}
Notice that since every $I_{f(\cdot ,s)}$ is an analytic set of pure codimension two such is 
their locally finite union $I_f$. These items force us to give the following definition.

\begin{nndefi} 
\label{mer-map2}
A meromorphic mapping $f: \calx \merto L\pp^1$ from a Hilbert manifold $\calx$ to the Sobolev 
loop space $W^{1,2}(\ss^1, \pp^1)$ is defined by a family of meromorphic on $\calx$ functions 
$\{f(\cdot ,s): s\in \ss^1\}$ such that: 

\sli the family of their indeterminacy sets $\{I_{f(\cdot , s)}:s\in \ss^1\}$ is locally
finite in $\calx$;

\slii $f$ is holomorphic when restricted to $\calx\setminus\bigcup_{s\in \ss^1} I_{f(\cdot, s)}$.
\end{nndefi}

Extension properties of meromorphic mappings with values in $L\pp^1$ occur to be not 
as good as for those with values in $\pp(l^2)$. In Example \ref{loop-exmp} we construct 
a holomorphic mapping $f: \cc^2\setminus \Gamma \to W^{1,2}(\ss^1,\pp^1)$, where $\Gamma = 
\{\gamma (s)= (\gamma_1(s),\gamma_2(s))\}$ is a loop in $\cc^2$ of class $W^{1,2}\subset 
\calc^{1,\frac{1}{2}}$, such that for every $s\in \ss^1$ the point $\gamma (s)$ is an 
essential singularity of $f$. That is $f$ doesn't extend to any neighbourhood of 
$\gamma (s)$ even meromorphically. This example shows that neither Thullen no Hartogs 
type extension theorems are valid for meromorphic  mappings with values in $L\pp^1$.

\begin{nnrema}  \rm
What concerns Theorem \ref{non-imb} assumption that there exists a meromorphic 
injection $g :L\pp^1 \merto \pp (\l^2)$ quickly leads to a contradiction between the 
already mentioned Example \ref{loop-exmp} and Corollary \ref{hart-corol}, see the 
end of the section \ref{LOOP} for more details.
\end{nnrema}

\newsect[MER-F]{Meromorphic functions on Hilbert manifolds}

\newprg[MER-F.h-m]{Hilbert manifolds and analytic sets}

For standard facts from infinite dimensional complex analysis we refer to \cite{Mu}. By a 
complex Hilbert manifold in this paper we shall understand a Hausdorff topological space 
$\calx$ covered by coordinate charts $(U_k,h_k)$, where $h_k:U_k\to V_k$ are 
homeomorphisms to open subsets of a separable complex Hilbert space $L$, \ie $L$ is 
$\cc^n$ or $\l^2$, such that the transition maps are holomorphic. We shall say, if needed, 
that $\calx$ is modelled over $L$. If $L=\cc^n$ we say that $\calx$ is $n$-dimensional, if 
$L=\l^2$ that $\calx $ is infinite dimensional.

\smallskip  We shall consecutively use results about Hilbert analytic sets from \cite{Ra}, 
where they are proved for analytic sets in Banach spaces in which every closed subspace admits
a closed complement. By $_L\calo_{z_0}$ we denote the ring of germs of holomorphic functions 
at $z_0$ in the separable complex Hilbert space $L$. For $L=\cc^n$ we write $_n\calo_{z_0}$, 
for $L=\l^2$ we write $_{\l^2}\calo_{z_0}$. But more often when  $L$ is clear from the context
we shall write simply $\calo_{z_0}$. This is a factorial domain, see Theorem I.1.4.4 in 
\cite{Ra}. Let $\calx$ be a Hilbert manifold.

\begin{defi}
\label{has-def}
By a Hilbert analytic set in $\calx$ we understand a subset $A\subset \calx$ such that 
for every point $a\in \calx$ there exists neighbourhood $\calv \ni a$ and a holomorphic
mapping $h:\calv \to F$ to a separable complex Hilbert space $F$ such that 
\[
A \cap \calv = \{z\in \calv : h(z) =0\}.
\]
\end{defi}
If $F=\cc^n$ for all $a\in A$ we say that $A$ is an analytic set of {\slsf finite definition}.
If $A$ is locally contained in a finite dimensional locally closed submanifold $M$ of $\calx$ and
is an analytic set in $M$ we say that $A$ is {\slsf finite dimensional}.

\smallskip
We denote the germ of an analytic subset $A$ at $a$ by $A_a$ and by $\calj(A_a) \subset \calo_a$ 
the ideal of germs of holomorphic functions that vanish identically on $A_a$. We say that $A$ is 
principal at $a$ if the ideal $\calj(A_a)$ is principal and we say that $A$ is principal if $A_a$ 
is principal at every $a \in A$. Given a holomorphic map $f: \calx \to F$ we write 
$V(f)= f^{-1}(0)$ and if $f=(f_1, \dots, f_n)$ is $\cc^n$-valued we write $V(f_1,\dots,f_n)= 
V(f_1) \cap \dots \cap V(f_n)$. Given $a \in \calx$, the zero locus of a holomorphic germ 
$V(f_a)$ is defined 
by the germ of analytic set at $a$ of the zero locus of $f$, i.e $V(f_a)=\big(V(f)\big)_a$.  
We say that a germ of an analytic subset $A_a$ is of pure codimension $p$ at $a$ if the height 
of the ideal $\calj(A_a)$ is equal to $p$. This is equivalent to saying that 
\[ 
\sup \{ n \in \mathbb{N} \ | \  \exists E \text{ subspace of } L , \  \dim E=n \text{ and }
\left(a + E \right) \cap A_a = \{a\}  \} =p, 
\] 
see corollary on page 70 of \cite{Ra}. And we say that $A$ is of pure codimension $p$ if 
$\codim(A_a)=p$ for all $a \in A$. It is proved in \cite{Ra}, see Proposition II.3.3.3 that 
for any $d \in \nn$ the germ $A_a$ of a Hilbert analytic set $A$ at $a$ admits the following 
decomposition 
\begin{equation}
\eqqno(decomp1)
A_a = \left( \bigcup_{k \leqslant d} A^k_a \right) \cup A_{a,d},
\end{equation}
where $A_a^k$ are irreducible germs of pure codimension $k$ for $k=1,\dots,d$ and $A_{a,d}$ 
is the union of germs of codimension $>d$. 

\smallskip For a Hilbert analytic set in $\calx$ the regular part $A^\text{reg}$ of $A$ is the 
set of points $a \in A$ such that $A$ is a submanifold of $\calx$ in a neighbourhood of
$a$. $A^\text{sing}:= A \setminus A^\text{reg}$ is called the singular locus of $A$. 
Clearly $A^\text{reg}$ is a union of arcwise connected components, which are locally closed 
submanifolds of $\calx$. Let $B^k$ be such a component of codimension $k$ (suppose 
such exists).  It is proved in \cite{Ra}, see Theorem II.4.1.2 that the closure $\bar B^k$ 
is a Hilbert analytic subset of $\calx$ of codimension $k$ which is globally irreducible in 
$\calx$. We call $\bar B^k$ an irreducible $k$-codimensional component of $A$. The union of 
all $\bar B^k$-s  is locally finite in $\calx$, see Proposition II.3.3.2 in \cite{Ra}.

\smallskip We shall use under several occasions the Weierstrass preparation theorem in
the Hilbert settings. The proof is the same as in the finite dimensional case. Let $f$ 
be a holomorphic function in open neighbourhood $\mathcal{U}$ of $a \in l^2$ such that 
$f(a)=0$ but $f \not\equiv 0$. Then there exists a direction $v \in l^2$ such that
the map $\lambda \mapsto f(a+\lambda v)$ does not vanish identically. Take a hyperspace
$E'$ in $l^2$ such that $\l^2 = E'\oplus \cc v$. 

\begin{thm}
\label{weier-th}
There exists a neighbourhood $\mathcal{V}$ of the origin such that for all $z=z'+\lambda v 
\in \mathcal{V}\cap (E' \oplus \mathbb{C}v)$ one has
\begin{equation}
\eqqno(weier1)
f(z+a)= \left( \lambda^k + \sum_{j=1}^{k }c_j(z')\lambda^{k-j} \right) \phi(z)
\end{equation}
with $\phi$ holomorphic and non-vanishing on $\mathcal{V}$. Moreover, for all $j \in 
\{1,k]\}$ the holomorphic coefficient $c_j \in \mathcal{O}(E')$ satisfies $c_j(0')=0$. 
\end{thm}
\proof 	By taking $g(z) = f(z+a)$ one can consider only the case where $a=0$ and by continuity of $g$ one can found a direction $v \in l^2 \setminus \{0\}$, a ball $B' \subset E'$ centred at the origin and a small disk $\Delta $ such that for all $(z', \lambda) \in \overline{B'} \times \overline{\Delta}$ one has $z'+\lambda v \in \mathcal{U}$ and for all $(z', \lambda) \in B' 
\times \partial \Delta$, $f(z'+\lambda v) \neq 0$. 	The order of the zero is given then by 
$k = \frac{1}{2\pi i } \int_{\partial \Delta} 	\frac{ \frac{\partial g}{\partial \lambda}(\xi v) }
{g(\xi v)}\, d\xi $

\smallskip Let us denote the zeros of $g$ by $\lambda_1(z'), \dots, \lambda_k(z')$. Then we 
obtain Newton polynomials evaluated on the zeros by
\[
	\sum_{j=1}^{k} \lambda_j(z')^\mu =
	\frac{1}{2\pi i } \int_{\partial \Delta} \xi^\mu
	\frac{ \frac{\partial g}{\partial \lambda}(z'+ \xi v) }{g(z'+ \xi v)}\, d\xi 
\]
From algebra one knows that elementary symmetric  polynomials can be expressed as polynomials
from Newton polynomials. Therefore the polynomial 
\[
	P(z'+\lambda v) = \prod_{j=1}^k (\lambda - \lambda_j(z')) =  \lambda^k + \sum_{j=1}^{k }c_j(z')\lambda^{k-j} 
\]
has holomorphic  in $z'$ coefficients. Now $\phi$ is defined by $\phi(z) = \frac{g(z)}{P(z)}$. 

\smallskip\qed 

\smallskip 
\begin{rema} \rm 
The monic polynomial $P(z) \in \mathcal{O}(E')[\lambda]$ given by $P(z) = 
\lambda^k + \sum_{j=1}^{k }c_j(z')\lambda^{k-j}$ is called a Weierstrass polynomial.     
\end{rema}

\smallskip If $A$ is an analytic set in a complex Hilbert manifold $\calx$ of 
codimension $\ge 2$ then $A$ doesn't contain a principal germ at any of its points.
This means in other words that in decomposition \eqqref(decomp1) $A_a^1=\emptyset$ for every 
$a\in A$. The following statement gives a more detailed description of analytic sets of 
codimension $\ge 2$. It is proved in \cite{Ra}, see Lemma II.1.1.12. 
We shall give a
proof of it as well as of the subsequent Proposition \ref{clos} in order to make our 
exposition more accessible.

\begin{lem} 
\label{codim2lem}
Let $A$ be a germ at zero of an analytic set in a complex Hilbert space $L$. If $A$ 
does not contain a principal germ at zero then there exist holomorphic germs  $f_1$ and $f_2$ 
at zero such that $A \subset V(f_1,f_2)$. Moreover $V(f_1,f_2)$ does not contain a principal 
germ at zero.
\end{lem}
\proof Since $A$ is analytic there exists $h \in \calo(L,F)$, with $F$ a complex Hilbert 
space, such that $A =V(h)$. Since $h\not\equiv 0$ we find a non-zero element $u_1$ of the 
dual $F^*$ such that $u_1 \circ h \not\equiv  0_{L^*}$. Let $u_1\circ h = g_1^{\alpha_1} 
\dots g_n^{\alpha_n}$ be the decomposition to irreducible factors and set $f_1= g_1 \dots g_m$. 
Choose a direction $\vect \in L$ such that $f_1$ does not vanish identically on $\<\vect \>$ 
and consider the decomposition $L= \<\vect\>^{\perp}\oplus\<\vect\>  $ with corresponding
coordinates $(z',z'')$. The Weierstrass Preparation Theorem applied to $f_1$ gives that $f_1= 
A P_1$ where $P_1(z',z'') \in \calo(\<\vect\>)[z'']$  is a Weierstrass polynomial and $A$ is
an invertible germ. Divide $h$ by $P_1$ using the Weierstrass Division theorem and get $h = 
P_1 Q +R$ where $Q \in 
\calo(L,F)$ and $R \in \calo(\<\vect\>^{\perp},F)[z'']$ with $\deg R < \deg P_1 $. 

\smallskip Now assume that for $u \in F^*$ we have $u\circ h |_{V(f_1)} = 0$. Then there exists 
$q \in \calo(L)$ such that $u \circ h = P_1 q$. At the same time we have $u\circ h = P_1 
\ (u \circ Q) + u \circ R$ and  we can deduce from the unicity of the Euclidean division that 
$u \circ R = 0$. If this holds for every $u\in F^*$ we obtain that $R=0$ and $h= P_1 Q$. Therefore 
$A$ would contain the principal germ $V(P_1)=V(f_1)$. Contradiction. We conclude that exists is 
$u_2 \in F^*$ such that $f_2 = u_2 \circ h$ satisfies $f_2|_{V(f_1)} \neq 0$. Now $A\subset V(f_1,f_2)$ 
and $V(f_1,f_2)$ doesn't contain a principal germ at zero as stated. 

\smallskip\qed

The following proposition is a particular case of Proposition III.2.1.3 from \cite{Ra}. 
\begin{prop}
\label{clos}
Let $\calx$ be a complex Hilbert manifold and $A$ an analytic subset of codimension $\geqslant 2$ 
of $\calx$. Let $C$ be an analytic set of pure codimension one in $\calx \setminus A$. Then 
$\overline{C}$ is an analytic subset of $\calx$.
\end{prop}
\proof The question is local, therefore the only thing we need to do is to prove that $\bar C$ 
is analytic in a neighbourhood of any point $a \in \overline{C} \cap A$. In the sequel we 
assume that $a=0$  
and let $\calu$ be a neighbourhood of $0$ in $L$. By Lemma \ref{codim2lem} we can suppose that 
$A = V(f_1,f_2)$. If for every direction $\vect\in L$ the  analytic  in $ \< \vect \> \cap 
(\calu \setminus A)$ set $\<\vect\> \cap \calu \cap C$ is not discrete in  $\<\vect\> \cap 
\left( \calu \setminus A \right)$ we have that $C$ contains an open subset of $L$, which 
is a contradiction. 

\smallskip  Therefore there exists a direction $\vect \in L$ such that $A \cap  \<\vect\> = 
\{ 0\}$ and $\<\vect \> \cap C $ is discrete  in $\calu \setminus A$. Decompose $L= L' 
\oplus \<\vect \>$ where $L' = \<\vect \>^{\perp}$. Taking $\calu$ smaller, if necessary, 
we can assume that $\calu= \calu' \times \calu''$ corresponding to this decomposition. Moreover,
we can chose $\calu'$ and $\calu^{''}$ in such a way that $\bar C\cap (\calu'\times \d\calu^{''})
=\emptyset$. Let $\pi: \calu \to \calu'$ be the orthogonal projection to $\calu'$ and $P_1$ and 
$P_2$ be Weierstrass polynomials of $f_1$ and $f_2$ respectively.  Denote by $\rho(P_1, P_2) \in 
\calo(U')$ the resultant of $P_1$ and $P_2$. We have that $\pi(A)= \{x \in \calu' \ | \ 
\rho(P_1,P_2) =0 \}$ and the latter is an analytic subset of $\calu'$ of codimension one. Set 
$\calu_0 = ( \calu' \backslash \pi(A)) \times \calu''$. Then $( \overline{A} \cap \calu_0 , 
\pi|_{\overline{C} \cap \calu_0}, \calu' \backslash \pi(A))$ is a ramified cover. 
Denote by $d$ the number of preimages of a generic $z' \in U' \backslash \pi(A)$. The fibre 
over $z'$ writes as
 
\[
 \pi^{-1}(\{z'\}) = \{z'+z_1''(z')\vect, \  \dots \ , \ z'+z_d''(z')\vect\}.
\] 
Consider the function $P: L' \times \cc \to \cc$ defined by

\begin{equation}
  P(x',z)= \prod_{j=1}^d \big(z- x_j''(x') \big).
\end{equation}
We have that $C\cap \calu_0 =\{ z'+z\vect \in \calu\backslash A \ | \ P(z',z)= 0 \} $.
Elementary symmetric functions of $\{z_j''(z')\}_{j=0}^d$ are bounded and holomorphic and 
therefore extend holomorphically to $\calu'$. Therefore $\bar C\cap \calu = 
\{ z'+z\vect \in \calu \ | \ P(z',z)= 0 \} $, \ie $\bar C$ is analytic in $\calu$
as stated. 

\smallskip \qed 

\newprg[MER-F.hol]{Holomorphic functions and mappings}
An $\l^2$-valued holomorphic mapping from an open subset $\calx$ of $\l^2$ is a mapping 
$f:\calx \to \l^2$ which is Fr\'echet differentiable at all points of $\calx$. It is 
well known that $f:\calx \to \l^2$ is holomorphic if and only if it is $G$-holomorphic, 
\ie G\^ateaux differentiable, and $\norm{f(z)}$ is locally bounded, see Proposition 8.6
and Theorem 8.7 in \cite{Mu}. This implies that if $f$ is written in coordinates as 
\begin{equation}
\eqqno(coord-hol1)
f(z) = (f_1(z), f_2(z),...)
\end{equation}
then $f$ is holomorphic if and only if all $f_k$ are holomorphic functions 
and $\norm{f(z)}^2=\sum_{k=1}^\infty \left| f_k(z) \right|^2$ is locally bounded. Indeed, 
since we need to prove the G\^ateaux differentiability
only we can suppose that $z\in \cc$ and then Cauchy formula and local boundedness of $\norm{f(z)}$ 
gave us the local boundedness of $\norm{f'(z)}$ where $f'(z)=(f_1'(z), f_2'(z), \dots)$. Indeed 
\begin{equation}
 \eqqno(coord-hol2)
\norm{f'(z)}^2 =  \frac{1}{(2\pi )^2}\sum_k\left|\int\limits_{|\zeta -z|=\eps}
\frac{f_k(\zeta)}{(\zeta -z)^2}d\zeta \right|^2 \le \frac{1}{(2\pi )^2\eps^4}\sum_k
\left(\int\limits_0^{2\pi }\left| f_k(\zeta)\right|\eps d\theta \right)^2 \le 
\end{equation}

\[
\le \frac{1}{(2\pi \eps )^2}\sum_k \int\limits_0^{2\pi }\left| f_k(\zeta)\right|^2 d\theta 
\cdot 2\pi \le
\frac{1}{2\pi\eps^2}\sup_{|\zeta - z|=\eps}\norm{f(\zeta)}^2.
\]
In the same fashion one proves the local boundedness of $\norm{f^{''}(z)}$ and hence the continuity 
of $f'(z)$. And this in its turn implies G\^ateaux differentiability:
\[
\norm{f(z+h)-f(z) -f'(z)h} = \norm{\int\limits_0^1\frac{d f(z+th)}{dt}dt-f'(z)h} = 
\]
\[
= \norm{\int\limits_0^1(f'(z+th)-f'(z))dt\cdot h} = o(h)
\]
by the continuity of $f'(z)$.

\begin{lem} 
\label{hol-ext}
Let $A$ be an analytic subset of the ball $B\subset L$ that doesn't contain a principal germ at any of its 
points. Then every holomorphic map $h \in \calo(B\setminus A ,F)$, where $F$ is a Hilbert space, extends 
holomorphically to $B$.   
\end{lem}

\proof It is sufficient to prove that $h$ extends to a neighbourhood of a given point, say a zero point.   
By the Lemma \ref{codim2lem} there exist germs $f_1$ and $f_2$ such that the germ of $A$ at zero is 
contained in $V(f_1,f_2)$. Choose a direction $a \in L$ such that $f_1$ and $f_2$ do not vanish 
identically on $\text{Vect}(a)$ and decompose $L=L' \oplus \text{Vect}(a)$.  The Weierstrass Preparation 
Theorem gives the Weierstrass polynomials $P_1(z',z)$ and $P_2(z',z)$ in $\calo(L')[z]$ such that $X 
\subset V(f_1,f_2)=V(P_1,P_2)$. Therefore there exists $r_1>0$ such that $\{0_{L'}\}\times \bar 
\Delta_{r_1} \subset B$ and $P_1(0_{L'},z)$ and $P_2(0_{L'},z)$ do not vanish on $\{0_{L'}\} \times 
\partial \Delta_{r_1}$. 
  
\smallskip The  number of roots of $P_1(0_{L'}, z)$ and $P_2(0_{L'}, z)$ is finite. By continuity of 
$P_1$ and $P_2$ there exists $r_2$ such that $P_1$ and $P_2$ do not vanish on $B(0_{L'},r_2) \times 
\partial \Delta_{r_1}$, i.e. $V(P_1, P_2) \subset B(0_{L'},r_2) \times \bar\Delta_{r_1}$.
We define $\tilde h$ by
\begin{equation} \label{int}
\tilde h(z',z) = \frac{1}{2 \pi i } \int_{\partial \Delta_{r_1}}h(z',u) \frac{du}{u-z}.
\end{equation}
This $\tilde h$  a holomorphic function on $B(0_{L'},r_2) \times \Delta_{r_1} $. Let 
$Q_1, ..., Q_m$ the irreducible factors of $P_1$ and consider the resultant $\rho_j = 
\text{Res}(Q_j, P_2)$. $Q_j$ does not divide $P_2$ because in the other case $V(Q_j)$ would be a 
principal germ contained in $V(P_1,P_2)$. Therefore $\rho_j \not\equiv 0$. Consider $Y= V(\rho_1 
\dots \rho_m)$. Let $\pi:L \to L'$ be the canonical projection. Since $\pi(A)$ is contained in $Y$,
$\tilde h$ coincides with $h$ on $B(0_{L'},r_2) \setminus \left(Y \times \bar \Delta_{r_1}\right)$. 
Therefore $\tilde{h}$ is the holomorphic extension of $h$ to a neighbourhood of $0$. 
 
\smallskip\qed 
 
\begin{rema} \rm
Using (\ref{int}) and the Riemann extension theorem one can prove that if $A$ is a proper analytic 
subset of a Hilbert manifold $\calu$ and $h$ is a bounded holomorphic function on $\calu \setminus A$ 
then $h$ extends holomorphically to $\calu$. 
\end{rema}

\newprg[MER-F.mer]{Meromorphic functions}
The field of quotients of $_L\calo_{z_0}$ is denoted as $_L\calm_{z_0}$  (or, simply as 
$\calm_{z_0}$) and is called the ring of germs of meromorphic functions at $z_0$. Every 
such germ is defined in an open neighbourhood $\calu_{z_0}$ of $z_0$ and will be called 
a {\slsf meromorphic fraction}. A meromorphic function on a Hilbert manifold $\calx$ is
given by a (sufficiently fine) open covering $\{\calu_k\}$ of $\calx$, where $\calu_k$ 
are connected, and meromorphic  fractions $f_k$ on $\calu_k$ such that $f_k = f_j$ on 
$\calu_k\cap \calu_j$ for all $k$ and $j$. Another open covering $\{\calv_k\}$ and 
meromorphic fractions $g_k$ define the same meromorphic function if $f_k = g_j$ on 
$\calu_k\cap \calv_j$. 

\smallskip If the covering $\{\calu_k\}$ is taken to be sufficiently fine we can find
holomorphic functions $h_k,g_k\in \calo (\calu_k)$ with $g_k\not\equiv 0$ 
such that $f_k = h_k/g_k$ and consequently $h_kg_j=h_jg_k$ on $\calu_k\cap \calu_j$.
By $\calm (\calx)$ we denote the field of meromorphic functions on  $\calx$. The 
following lemma is well known in finite dimensions and its proof in Hilbert settings 
is similar.

\begin{lem}
\label{rel-prime}
Let $f$ and $g$ two holomorphic functions on $a \in l^2$ and that are relatively prime 
on the ring $\mathcal{O}_a$ of holomorphic functions on a neighbourhood of $a$. Then 
there exists a neighbourhood $\mathcal{U}$ of $a$ such that $f,g \in \mathcal{O}(\mathcal{U})$ 
and are relatively prime on $\mathcal{O}_y$ for every $y\in \mathcal{U}$.   
\end{lem}
\proof By translation one can suppose $a=0$. 
Let \( v \in l^2 \setminus \{0\} \) such that $f$ and $g$ do not vanish identically on 
$\mathbb{C}v$ and consider the decomposition \( l^2 = E' \oplus \mathbb{C}v \). By 
multiplying \( f \) and \( g \) by invertible holomorphic functions in \( \mathcal{O}_0 \),
we may apply the Weierstrass Preparation Theorem in the direction \( v \). Then \( f \) and 
\( g \) can be written as Weierstrass polynomials in \( \mathcal{O}_{0'} (E')[\lambda] \).
Since \( f \) and \( g \) are relatively prime there exist \( u, v \in \mathcal{O}_{0'}(E')
[\lambda] \) such that \(	fu + gv = r \in \mathcal{O}_{0'}(E')\). 	This identity holds 
in a neighbourhood \( \mathcal{U} \) of \( 0 \), so now suppose that at some \( y \in 
\mathcal{U} \), the germs \( f_y \) and \( g_y \) have a common non-invertible  factor 
\( h \in \mathcal{O}_y \). Then \( h \) divides \( r \in \mathcal{O}_{y'}(E')[\lambda] \), 
where \( y = y' + \lambda_y v \). 
	
\smallskip Applying again the Weierstrass Preparation Theorem to \( h \) we can write 
\( h = u \cdot h_1 \), where \( u \in \mathcal{O}_y \) is invertible and \( h_1 \in 
\mathcal{O}_{y'}(E')[\lambda] \). Since \( h_1 \) divides \( f \), and \( f \) is monic 
in $\lambda$, it follows that $h_1$ must be invertible in \( \mathcal{O}_{y'}(E') \), 
hence \( h \) is invertible in \( \mathcal{O}_y \). Therefore \( f_y \) and \( g_y \) 
are relatively prime in \( \mathcal{O}_y \). This proves that there exists a neighbourhood 
\( \mathcal{U} \) of \( 0 \) such that for all \( y \in \mathcal{U} \), the germs 
\( f_y \) and \( g_y \) are relatively prime in \( \mathcal{O}_y \).

\smallskip\qed

\begin{rema} \rm
Due to this lemma we can refine a covering $\{\calu_k\}$ involved in the definition of a 
meromorphic function and assume in the sequel that $g_k$ and $f_k$ are relatively prime
at every point of $\calu_k$.
\end{rema}

The following lemma is usually attributed to Hadamard.
\begin{lem} 
\label{hadamard}
A formal power series with coefficients in the integral domain $\calk$
\[F(\lambda) = \sum_{n=0}^\infty a_{-n}\lambda^{-n}  \in \calk[[\lambda^{-1}]] \]
represents a rational function $\frac{P(\lambda)}{Q(\lambda)}$ with $P, Q \in 
\calk[\lambda]$ and 
$\deg Q \leqslant N$ if and only if 
\begin{equation}
\label{det}
\left|  \begin{array}{cccc} a_{-n_1} & a_{-n_2} & \dots & a_{-n_{N+1}} \\ \vdots & \vdots & 
\vdots & \vdots \\ a_{-n_1-N} & a_{-n_2 -N} & \dots & a_{-n_{N+1}-N}  \end{array} \right|
= 0  \end{equation} 
for all $(N+1)$-tuples $1\le n_1 < \dots < n_{N+1}$.
\end{lem}
\proof
Indeed we look for a non-zero polynomial $Q(\lambda)= c_0+c_1 \lambda + \dots +c_N\lambda^N$ 
with coefficients in $\calk$ such that $FQ \in \calk[\lambda]$. But this condition means that 
for every $k \geqslant 1$ one should have 
\begin{equation}
\label{hyp-plane}
a_{-k}c_0+ \dots + a_{-k-N}c_N=0.
\end{equation}
This relation means that vectors $b_k:=(a_{-k}, a_{-k-1}, \dots, a_{-k-N})$, 
$k \in \nn^*$, belong to the hyperplane with the equation \ref{hyp-plane} in the 
$\calm$-linear space $\calm^{N+1}$, here $\calm$ is the field of fractions of $\calk$. 
The latter means that every $N+1$ of them are linearly dependent, and this is precisely 
what tells the condition \ref{det} 

\smallskip\qed

Let $A$ be a subset of $B^\infty$. We say that $A$ is thick at $z_0 \in B^\infty$ if for any neighbourhood 
$\calv \ni z_0$ the set $A \cap \calv$ is not contained in a proper analytic subset of $\calv$. The domain $R^\infty_{r_1,r_2}$ 
\begin{equation}\label{ringdomain}
R^\infty_{r_1,r_2}=A_{r_1,r_2}\times B^\infty
\end{equation}
 will be called a ring domain. Here $A_{r_1,r_2}= \Delta_{r_1}\setminus \bar \Delta_{r_2}$ is an annulus,  $0 \leqslant r_1 < r_2$. The following lemma is an infinite dimensional version of  the celebrated Levi's theorem, \cite{Lv}.

\begin{lem} 
\label{levi-th-n}
Let $f$ be a holomorphic function in the ring domain $R_{1-r,1}^\infty$. Suppose that for 
$w$ in some subset $A\subset B^\infty$ thick at origin restrictions $f_w:=f(\cdot , w)$ 
meromorphically extend from $A_{1-r,1}$ to $\Delta$ and the number of poles counted with 
multiplicities of these extensions is uniformly bounded. Then $f$ extends to a meromorphic 
function on $\Delta \times B^\infty$. 
\end{lem}
\proof Write $f$ as 
\begin{equation}
f(z,w) = f^+(z,w) +f^-(z,w) = \sum_{n \geqslant 0 }a_n(w)z^n + \sum_{n<0}a_n(w)z^n
\end{equation}
where $a_n \in \calo(B^\infty)$. Notice that $f^+$ is already holomorphic in $\Delta \times 
B^\infty$. Our task therefore is to extend $f^-$. By Lemma \ref{hadamard} applied to the ring 
$\cc$ the extendability of $f^-(z,w)$ for $w \in B^\infty$ to the disk together with the 
condition on the poles means that for $a_n = a_n(w)$ the determinants \ref{det} vanish. 
Therefore they vanish identically as functions of $w$. And therefore, again by Lemma 
\ref{hadamard} but this time applied to the ring $\calo(B^\infty)$ , we see that $f^-(z,w)$
is rational over the field $\calm(B^\infty)$. I.e., is meromorphic on $\Delta \times 
B^\infty$.

\smallskip\qed

\begin{rema} \rm
\label{levi-open}
Lemma \ref{levi-th-n} can be restated without condition on boundedness of the number of
poles, this is convenient in applications. The assumption should be that $f(\cdot,w)$ 
meromorphically extends from $A_{1-r,1}$ to $\Delta$ for $w$ from the set $A$ which is 
not contained in a countable union of proper analytic sets, ex. $A$ is open and 
non-empty. 
\end{rema}

\begin{lem} 
\label{ext-mero-f}
Let $I$ be an analytic set with $\codim I \geqslant 2$ in a complex Hilbert manifold $\calx$ and let $f$ be a meromorphic function on $\calx \backslash I$. Then $f$ extends to a meromorphic function $\tilde{f}$ on the 
whole of $\calx$.
\end{lem}
\proof 
Denote $\calp_f$ the divisor of poles of $f$. $P_f \subset \calx \backslash I$ is a complex 
hypersurface and by Proposition \ref{clos} its closure of $\overline{P_f}$ is an analytic
subset of $\calx$. Let $a \in I$. If there is a neighbourhood $U$ of $a$ such that 
$\overline{P_f}\cap U = \emptyset$ then $f \in \calo(U \backslash I)$ extends holomorphically 
to $U$  by Lemma \ref{hol-ext}. Let $a \in I\cap \overline{P_f} $. Choose coordinates $(z,w)$ 
in a neighbourhood $U$ of $a=0$ associated to the decomposition $L = \<b\> \oplus L'$ 
for some $b \in L$ such that $\<b\>\cap \bar P_f$ is discrete. Now one can place a ring domain 
$R_{1-r,1}^\infty \subset U\setminus \bar P_f$, \ie $f$ will be holomorphic on $R_{1-r,1}^\infty$. 
Note that for every $w \in B^\infty \setminus \pi (I)$ function $f(\cdot,w)$ extends 
meromorphically to $\Delta$. Here $\pi : \<b\> \oplus L'\to L'$ is the natural projection. 
$B^{\infty}\setminus \pi (I)$ is not contained in a countable union of proper analytic sets,
otherwise $\Delta\times B^{\infty}$ would be such. By Lemma \ref{levi-th-n}, taking  into 
account Remark \ref{levi-open}, $f$ extends meromorphically to a neighbourhood of $0$. 
Therefore $f$ extends meromorphically to $\calx$.
 
\smallskip\qed

\smallskip
We continue with a Thullen type extension of meromorphic functions. 

\begin{thm} \label{thullen}
Let $A$ be a proper analytic subset of a complex Hilbert manifold $\calx$ and let $f$ be a meromorphic function on $\calx \setminus A$. Suppose that for every irreducible component 
$A_1$ of $A$ of codimension one there exists a point $p_1 \in A_1$ such that $f$ 
meromorphically extends to some neighbourhood $V_1$ of $p_1$. Then $f$ meromorphically 
extends onto the whole of $\calx$. 
\end{thm}
\proof
If $A$ has codimension $\geqslant 2$ irreducible components then by Lemma \ref{ext-mero-f} $f$ extends across them. Now consider $A_1$ an irreducible component with $\codim A_1 =1$. Let $p_1$ be a point where $f$ extends and take $p_2$ another point of $A_1$. Take another point 
$p_2$ which is not singular for $A_1$ and doesn't belong to the intersection of $A_1$ 
with other components of $A$. Consider a path $\gamma: [0,1] \to A_1$ such that $\gamma(0)=p_1$ and $\gamma(1)=p_2$. Choose $\gamma$ in $\Reg A_1$ minus intersections with other components
of $A$. Define the set $E$ by 
\begin{equation}
E= \{ t \in [0,1] \ | \ f \text{ extends meromorphically to a neighbourhood of } \gamma([0,t]) \}.
\end{equation}
Let $t_0=\sup E$. We need to prove that $t_0=1$. Suppose not. Then one can find a complex disk $\Delta$ through $\gamma(t_0)$ transverse to $A_1$ that does not intersect other components of $A$ and is not contained in the set of indeterminacy of $f$. 
We can find a local coordinates $(\lambda,w)$ with center $\gamma(t_0)$ such that $\Delta \times \{0 \}$ is our disk. By Lemma \ref{levi-th-n}  $f$ extends meromorphically to a neighbourhood of $\gamma(t_0)$. Contradiction. Therefore $t_0=1$ and $f$ extend meromorphically across $A_1$
minus a codiomension $\ge 2$ subset. Do the same for all other components. We get an 
extension of $f$ onto $\calx$ minus a codimension $\ge 2$ analytic subset. Now Lemma 
\ref{ext-mero-f} applies and finishes the proof.

\smallskip\qed

For the definition of the infinite dimensional $1$-concave Hartogs figure 
see \eqqref(inf-hart). 


\begin{cor}
\label{mer-hart}
Every $f \in \calm( H^\infty_1(r))$ extends meromorphically to $\Delta \times B^\infty$. 
\end{cor}
\proof  Notice first that the set $I$ of such $w \in B^\infty$ that $I_f \supset A_{1-r,1}
\times \{w\} $ is a proper analytic subset of $B^\infty$. One easily derives local equations 
for $I$ from that of $I_f$.  If $w_0 \in I$ then for $z_0 \in A_{1-r,1}$ the equation of $I_f$ 
in a neighbourhood of $(z_0,w_0)$ has the form $f(z,w)= \sum_{k=0}^\infty f_k(w)(z-z_0)^k=0 $, 
where $f_k$ are $\cc^2$-valued holomorphic functions in a neighbourhood of $w_0$.
Notice that all $f_k(w_0)=0$ and therefore in a neighbourhood of $w_0$ we have 
$I=\{f_1(w)= f_2(w)= \dots =0 \}$. Extension of $f$ 
to $\Delta \times \left( B^\infty \setminus I \right) $ can be achieved as in the proof of 
Theorem \ref{thullen} using Theorem \ref{levi-th-n}. Now $f$ is meromorphic on $(\Delta \times 
B^\infty) \setminus (\Delta \times I)$ and meromorphic on $A_{1-r,1} \times B^\infty$. So we 
are under assumptions of the Thullen-type theorem \ref{thullen} and conclude that $f$ extends 
to a meromorphic functions on $\Delta \times B^\infty$.  

\smallskip\qed

\newsect[PROJ]{Meromorphic mappings to the projective  Hilbert space}

\newprg[PROJ.hol]{Holomorphic mappings to the projective Hilbert space}

\smallskip In the definition of the projective Hilbert space over $\l^2$ one requires homogeneous
coordinates $[w]=[w_0:w_1: ...]$ to be such that $w=(w_0,w_1,...)\in \l^2\setminus \{0\}$. Notice that 
the affine coordinate $(z^j_1,z^j_2,...)$ in the chart $\calu_j$ belong to $\l^2$ as well. Denote by 
$h_j: \calu_j\to \l^2$ the corresponding coordinate map, namely $h_j([w])= \{\frac{w_k}{w_j}\}_{k \neq j}$. 
Transition maps 
\[
h_j\circ h_i^{-1} : \l^2\setminus \{z^i_j= 0\} \to \l^2\setminus \{z^j_i=0\} 
\]
are given by 
\begin{equation}
\eqqno(trans-map)
h_j\circ h_i^{-1} : (z^i_0,z^i_1,...) \mapsto (z^i_0/z^i_j,z^i_1/z^i_j,...) = (z^j_0,z^j_1,...).
\end{equation}
They are obviously biholomorphic and therefore $\pp (\l^2)$ is a complex Hilbert manifold 
modelled over $\l^2$. Closed submanifolds of $\pp (\l^2)$ will be called projective Hilbert 
or simply projective manifolds. I.e, a projective Hilbert  manifold is a closed subset 
$\calm$ of $\pp(l^2)$ such that for every point $m \in \calm$ there exists a neighbourhood 
$\calu \ni m$ and a holomorphic map $h: \calu \to L$ to a Hilbert space $L$ with surjective 
differential such that $\calm = \{ z \in \calu  :  h(z)=0 \}$. Their intersections with 
$\pp (l^2)\setminus \{w_0=0\}$ will be called affine. As well as closed 
submanifolds of $\pp (l^2)\setminus \{w_0=0\}$ will be called Stein.

\smallskip 
Let $\calx$ be a complex Hilbert manifold and let $f:\calx\to \pp (\l^2)$
be a holomorphic mapping. Without loss of generality we assume that the image of 
$f$ is not contained in any hyperplane $\calh_j$, see \eqqref(hyp-hj). Otherwise we would get 
$f:\calx\to \calh_j \cong \pp (\l^2)$ and so on. Taking into account the fact that
\[
\pp \deff \bigcap\limits_{\calh_j\supset f(\calx)} \calh_j
\]
is a closed linear subspace of $\pp (\l^2)$ isomorphic either to $\pp (\l^2)$ or to $\pp^N$
for some $N\in \nn$, we find ourselves under our assumption with $\pp$ as a target space of $f$. 
Set $A_j \deff f^{-1}(\calh_j)$. These are hypersurfaces (some of them may be empty) in $\calx$. 
For the restriction $f|_{\calx\setminus A_0}:\calx\setminus A_0 \to \calu_0$ the composition 
$h_0\circ f$ is a holomorphic $\l^2$-valued mapping. In affine coordinates $z^0_1,...,z^0_k,...$
it writes as
\begin{equation}
(h_0\circ f)(z) = (f^0_1(z),f^0_2(z),...),
\end{equation}
where $f^0_k$ are holomorphic and $\norm{(h_0\circ f)(z)}^2 = \sum_k|f^0_k(z)|^2$ 
is locally 
bounded in $\calx\setminus A_0$. 
\begin{prop}
\label{map1}
The functions $\{f_k^0\}$, holomorphic on $\calx\setminus A_0$, are in fact meromorphic on 
$\calx$. Moreover, their orders of poles along $A_0$ are locally bounded and consequently 
our mapping in a neighbourhood of any point $a\in \calx$  writes as 
\begin{equation}
\eqqno(map-phi1)
f(z) = [\phi_0(z):\phi_1(z):...],
\end{equation}
where $\phi_k$ are holomorphic, have no common zeros and, in addition, $\norm{\phi(z)}^2\deff 
\sum_k|\phi_k(z)|^2$ is locally bounded.
\end{prop}
\proof Take a point $a\in A_0$ such that $a\notin A_1$ \ie $f(a)\notin \calh_1$.
Then the holomorphic map $h_1\circ f : \calx\setminus A_1\to U_1$ writes as $h_1\circ f =
(f^1_0,f^1_2,...)$ with $f^1_0 = 1/f^0_1, f^1_2=f^0_2/f^0_1,..., f^1_k=f^0_k/f^0_1,... $  
holomorphic in a neighbourhood of $a$ and $\norm{(f^1_0(z),f^1_2(z),...)}$ locally bounded there. 
This proves that $f^0_k$ are meromorphic in a neighbourhood of $a$ and therefore on $\calx
\setminus (A_0\cap A_1)$. The same argument can be repeated for any $A_j$ and since 
$\bigcap A_j = \emptyset $ we see that $f^0_k$ are meromorphic on the whole of $\calx$.

\smallskip To prove the second assertion fix again a point $a\in A_0$ and let $h=0$ be the 
primitive equation of $A_0$ in a  neighbourhood $U$ of $a$. Existence of such $h$ follows from 
the fact that the ring $\calo_a$ is factorial. We can choose 
$U$ in such a way that $U = \Delta \times B$, where $B$ is the unit ball and $A_0\cap (\d \Delta 
\times B) = \emptyset$. Such $U$ will be called {\slsf adapted}
to $A_0$. 

\smallskip\noindent{\slsf Claim 1.} {\it For every $k\in \nn$ there exists $N_k\in \nn_0\deff\{0\}
\cup \nn$ such that 
\begin{equation} \label{not/h}
f^0_{k} = \frac{\psi_k}{h^{N_k}} \quad\text{ with } \psi_k \text{ holomorphic and not divisible by } 
h \text{ in } U.
\end{equation}
}

\noindent To prove this claim suppose first that $A_0\cap U$ is irreducible. Let $P(z_1,z')$ be 
the Weierstrass polynomial for $h$, \ie $h(z_1,z') = z_1^p+ a_1(z')z_1^{p-1}+...+a_p(z')$ with 
$a_j\in \calo (B), a_j(0')=0$ and $h = P\phi$, where $\phi$ is holomorphic in $U$ and doesn't 
vanish. Since $f^0_k$ has poles at most on $A_0\cap U$ we can multiply it by $P^{N_k}$, where 
$N_k\in \nn_0$ is the order of the pole of $f^0_k$ along $A_0\cap U$ to get a holomorphic 
function $\psi_k$ in $U$ which is obviously not divisible by $h$. If $A_0\cap U$ is reducible 
represent it as a union of a finite number of irreducible components $A_0\cap U= \bigcup_lA_0^l$ 
with primitive equations $h^l=0$ and repeat the argument as above to each $A_0^l$. The rest is 
obvious.

\smallskip\noindent{\slsf Claim 2.} {\it  Next we claim that the sequence $\{ N_k \}$ in 
\ref{not/h} is bounded.} If not one can find a subsequence $\{k_n\}$ such that $N_{k_n}\nearrow 
+\infty $ and 
\begin{equation}
f^0_{k_n} = \frac{\psi_{k_n}}{h^{N_{k_n}}} \quad\text{ with } \psi_{k_n} \text{ holomorphic an not divisible by } 
h \text{ in } U.
\end{equation}
Since $\bigcap A_j=\emptyset$ there exists $m$ such that $a\in A_0\setminus A_m$. Then 
\[
f^m_j = \frac{f^0_j}{f^0_m} = \frac{\psi_j}{\psi_m}\cdot h^{N_m-N_j}
\]
should be holomorphic in a neighbourhood of $a$ for all $j$. But this is certainly not true
for $j=k_n$ with $n$ big enough, contradiction. 

\smallskip Let $N= \max_{k \in \nn}{N_k}$ and set $\phi_0\deff h^N$, $\phi_k \deff f^0_kh^N$ for 
$k=1,2,...$. Consider the map
\begin{equation}
\eqqno(map-phi2)
\phi \deff (\phi_0,\phi_1,...) : U\to \l^2.
\end{equation}
We need to prove that $\phi$ indeed takes its values in $\l^2$ and is locally bounded on $U$. 
Remark that $\norm{\phi (z)} = |h(z)|^N\norm{\left(1,(h_0\circ f)(z)\right)}$ for $z$ in a neighbourhood 
of $\d \Delta \times B$ and therefore is {\slsf locally} bounded there. Using compactness of 
$\d \Delta$ we deduce that $\norm{\phi}$ is bounded on $\d\Delta \times B_{\delta}$, say by $M$, for 
$\delta >0$ small enough. Each cut $\sum_{k=0}^N|\phi_k|^2$ is plurisubharmonic on $U$ and 
therefore is bounded by $M$ on $\Delta\times B_{\delta }$. Tending $N\to \infty$ we get the boundedness 
of $\norm{\phi}$ on $\Delta\times B_{\delta}$. 

\smallskip All what left to prove is that $\phi_0,\phi_1,...$ can be taken without common zeroes.
Set $I=\{z\in U:\phi_0(z)=\phi_1(z)=...=0\}$. If $\phi_j$-s have common divisor in $_L\calo_a$ we can 
divide them by this divisor and this will not change our map $f$. Therefore $\codim I\ge 2$. Suppose
it is non-empty an choose some $b\in I$. Our map on $U\setminus I$ is given by 
\[
f(z) = [\phi_0(z):\phi_1(z):...]. 
\]
Since $f$ is holomorphic in a neighbourhood of $b$ there exists $j$ such that $f(b)\in \calu_j$. 
Therefore $\phi_j(b)\not=0$. Contradiction.

\smallskip\qed

\begin{rema} \rm 
\label{map2}
1. The affine coordinate chart $\calu_0$ will be considered as distinguished and $f^0_k$-s will 
be denoted simply as $f_k$-s. We proved that our holomorphic map $f$ is represented by a sequence 
of meromorphic on $\calx$ functions $(f_1,f_2,...,f_k,..)$ which are holomorphic on $\calx\setminus 
A_0 = f^{-1}(\calu_0)$ and $\norm{f}^2 = \sum_{k=1}^\infty{\norm{f_k}}^2$ is locally bounded 
on $\calx\setminus A_0$.

\smallskip\noindent 2. Moreover, we proved that for every point $a\in \calx$ there exits a neighbourhood 
$\calv\ni a$, holomorphic in $\calv$ functions $\phi_0,\phi_1,...$ without common zeroes such that
$\phi=(\phi_0,\phi_1,...)$ is locally bounded as a mapping $\calv \to \l^2$ and such that our $f$ is 
represented on $\calv$ as in \eqqref(map-phi1), \ie $f = \pi \circ \phi$, where $\pi :\l^2\setminus \{0\}
\to \pp (\l^2)$ is the canonical projection. 

\smallskip\noindent 3. Remark that if a holomorphic mapping $f:\calx\to \pp (\l^2)$ is given on 
an open $\calv\subset \calx$ by a representation \eqqref(map-phi1) then $A_j\cap \calv = \{\phi_j=0\}$. 
Indeed, $A_j\cap \calv$ is the preimage of $\calh_j$ under $\phi$, here $\calh_j = \{w\in \l^2: w_j=0\}$, 
since it is the preimage of $\calh_j$ under $f$.
\end{rema}

\newprg[PROJ.fin]{Meromorphic mappings to $\pp (\l^2)$.}

Recall that a meromorphic mapping between finite dimensional complex manifolds $X$ and $Y$ is a 
holomorphic mapping $f:X\setminus I \to Y$ where $I$ is analytic set with $\codim I \geqslant 2$, 
such that the closure $\bar\Gamma_f$ of its graph $\Gamma_f$ is an analytic set in $X\times Y$
and the restriction of the natural projection $\pi_1:X\times Y\to X$ to $\bar\Gamma_f$ is proper.
With some ambiguity this meromorphic mapping is denoted still as $f$, instead of $\bar\Gamma_f$ 
one writes simply $\Gamma_f$ and calls it the graph of $f$.

\smallskip Therefore an object which pretends to be a meromorphic map in Hilbert case 
should be at least holomorphic outside of an analytic set of codimension two.

\begin{prop}
\label{hol-ext1}
Let $\calx$ be a complex Hilbert manifold, $I\subset \calx$ an analytic subset 
of codimension $\ge 2$ and let $f:\calx \setminus I\to \pp (\l^2)$ be a holomorphic mapping. 
Then for every $a\in I$ there exists a neighbourhood $a\in \calu\subset \calx$ and 
functions $\{\phi_k\}$ holomorphic in $\calu$ such that 

\sli $\phi = (\phi_0,\phi_1, ...)$ is a holomorphic mapping from $\calu$ to $\l^2$;

\slii $\calu\cap I \supset \{z\in \calu: \phi_k(z)=0 \text{ for all } k\}$;

\sliii $(\pi \circ \phi)(z) = f(z)$ for $z\in \calu\setminus I$, where $\pi :\l^2\setminus \{0\}
\to \pp (\l^2)$ is the canonical projection.
\end{prop}
\proof As we saw in Proposition \ref{map1} mapping $f$ is defined by a sequence 
$\{f_k \}_{k \in \nn}$ of meromorphic functions on $\calx \backslash I$ which are holomorphic 
on $\calx \setminus (I \cup A_0)$ where $A_0=f^{-1}(\calh_0)$. We can extend $f_k$ meromorphically 
to $\calx$ by Lemma \ref{ext-mero-f}, denote by $\tilde f_k$ these extensions, and by 
Proposition \ref{clos} we can extend $A_0$ to an analytic set $\tilde A_0$ in $\calx$. Given 
$a \in I$ we need to examine $f$ in a neighbourhood of $a$. If $a\in I\setminus \tilde A_0$ 
all $f_k$ extend holomorphically to a neighbourhood of $a$ by Lemma \ref{hol-ext}, 
\ie our mapping $f$ extends holomorphically to a neighbourhood of such $a$. Suppose now that 
$a\in I\cap \tilde A_0$. Meromorphic functions $\tilde f_k$ have poles at most on $\tilde A_0$ 
and we need to prove that their orders stay bounded, \ie that there exists $N_k$ such that $h^{N_k} 
\tilde{f}_k$ is holomorphic on a neighbourhood of $a$ for every $k \in \nn$ and $N=\sup_{k \in \nn} 
N_k <\infty$. Here $h$ is a minimal defining function of $\tilde A_0$ at $a$. We assume that
that $N_k$ are minimal that do the job. Suppose this is not true, \ie there exists a subsequence
$N_{k_j}\to \infty$ such that $h^{N_{k_j}}\tilde{f}_{k_j}$ are holomorphic but 
$h^{N_{k_j}-1}\tilde{f}_{k_j}$ are not. But $h=0$ is an equation for $\tilde A_0$ in a 
neighbourhood of $a$ and we get in this case that $h^{N_{k_j}-1}\tilde{f}_{k_j}$ are not
holomorphic in this neighbourhood. This contradicts to Proposition \ref{map1}.
Now set $\phi_0= h^N$ and $\phi_k= h^N \tilde{f}_k$. Then \slii and \sliii are satisfied. The 
same consideration as in Proposition \ref{map1} gives us \sli. 

\smallskip\qed

The discussion above leads us to the following

\begin{defi}
\label{def-mer-p}
Let $\calx$ be a complex Hilbert manifold. A meromorphic map $f:\calx \merto \mathbb{P}(l^2)$  
is given by a holomorphic map $f: \calx \backslash I \to \mathbb{P}(l^2)$, where 
$I\subset \calx$  in an analytic subset of codimension $\geqslant 2$, such that for every 
$a \in I$ there exists a neighbourhood $\calv$ of $a$ and a 
holomorphic map $\phi =\{\phi_k\}_{k \in \nn}: \calv \to l^2$  such that $f= \pi \circ \phi$, 
\ie $f(z)= [\phi_0(z) : \phi_1(z) : \dots ]$ for $z \in \calv \backslash I$.  
\end{defi}
\noindent  The minimal such $I$ is called the 
indeterminacy set of $f$ and is denoted as $I_f$. We can interpret Proposition \ref{hol-ext1} 
as follows : if $f: \calx \setminus I \to \pp(l^2)$ is a holomorphic map and $\codim I 
\geqslant 2$ then $f$ extends meromorphically to $\calx$. One more remark: every meromorphic 
map $f: \calx \merto \pp(l^2)$ can be represented by globally defined meromorphic functions 
$\{f_k\}$ such that the orders of their denominators are locally uniformly bounded in the 
following sense. Let $A_0 = f^{-1}(\calh_0)$ be the preimage of the hyperplane $\calh_0 = 
\{z_0=0\}$ under the holomorphic map $f:\calx\setminus I_f\to \pp(\l^2)$. Then $A_0$ extends 
to an analytic hypersurface in the whole of $\calx$, we denote it still as $A_0$. And now the
poles of all $f_k$ are on $A_0$ and their orders are uniformly bounded in a neighbourhood 
of every point of $a$.

\smallskip
The indeterminacy set $I$ of a meromorphic mapping between Hilbert manifolds might look not as 
an analytic set in finite dimensions. And, moreover, the graph of $f$ might be not a direct 
analytic subset of $\calx \times \caly$ in the sense that it cannot be locally represented as 
a finite analytic cover. 

\newprg[PROJ.exmp]{Example}
Consider the following mapping $f:\l^2\merto \pp (\l^2)$
\begin{equation}
\eqqno(indet1)
f : z=\{z_k\}_{k=1}^{\infty} \mapsto \left[\{z_k(z_k-1/k)\}_{k=1}^{\infty}\right],
\end{equation}
or, in more details 
\[
\{z_k\}_{k=1}^{\infty} \mapsto  \left[z_1(z_1-1):...: z_k(z_k-1/k):...\right].
\]	
The indeterminacy set $I$ of $f$ (\ie all $z$ such that all components in \eqqref(indet1) 
vanish) is of Cantor type, see pp. 33-34 in \cite{Ra}. In particular it is totally 
disconnected and contains the sequence $\{e_1,...,\frac{1}{k}e_k,...\}$ converging to zero.

\smallskip 
Fix $z=\{z_j\}_{j=1}^{\infty}\in l^2$ and fix $k \in \nn^*$. Consider $y_{n,k} \in l^2$ defined 
for $n\in \nn^*$ by
\begin{eqnarray*}
y_{n,k} &=& \left( \frac{1}{k} + \frac{k z_k}{n^2} \right) e_k + \sum_{j=1, j \neq k}^n 
\frac{-jz_j}{n^2}e_j = \\
&=&\left( \frac{-z_1}{n^2}, \frac{-2z_2}{n^2}, \dots , \frac{-(k-1)z_{k-1}}{n^2},\frac{1}{k}+ \frac{kz_{k}}{n^2},  \frac{-(k+1)z_{k+1}}{n^2}, \dots,  \frac{-nz_{n}}{n^2}, 0 , \dots  \right).
\end{eqnarray*}
Notice that $y_{n,k} \underset{n \to \infty}{\longrightarrow} \frac{1}{k}e_k$. 
For $n$ big enough we have
\begin{eqnarray}
f(y_{n,k})&=& \left[\frac{-z_1}{n^2}\left(\frac{-z_1}{n^2} -1\right) : \frac{-2z_2}{n^2}\left(\frac{-2z_2}{n^2} - \frac{1}{2}\right): \dots :  \left( \frac{1}{k} + \frac{k z_k}{n^2} \right) \frac{kz_k}{n^2} : \dots \right] \\
&=&  \left[-z_1\left(\frac{-z_1}{n^2} -1\right) : -2z_2\left(\frac{-2z_2}{n^2} - \frac{1}{2}\right): \dots :  \left( \frac{1}{k} + \frac{k z_k}{n^2} \right) kz_k : \dots \right].
\end{eqnarray} 

\noindent
Therefore $f(y_{n,k})\underset{n \to \infty}{\longrightarrow}\left[z_1 : z_2 : \dots : z_k : \dots \right]$ and consequently $(y_{n,k}, f(y_{n,k})) \underset{n \to \infty}{ \longrightarrow} 
(\frac{1}{k}e_k, [z])$. This means that $\pi_1^{-1}\left(\left\lbrace\frac{1}{k}e_k \right\rbrace \right) = \pp(l^2)$ for every $k \in \nn^*$, here $\pi_1$ is the restriction of the
natural projection $l^2\times \pp (l^2)\to l^2$ to the graph of $f$. We conclude that 
the set $\pi_1^{-1}(I) \supset \{\frac{1}{k}e_k\} \times \pp(l^2)$ is not direct. 

\newprg[PROJ.rest]{Restrictions to complex curves}

Let us prove that any meromorphic mapping $f : \calx \merto \pp(l^2)$ restricted to a complex 
curve, which is not contained in the indeterminacy set of $f$, extends holomorphically to this 
curve. Since $f$ is meromorphic on $\calx$ there exists an analytic subset $I\subset \calx$ of 
codimension $\geqslant 2$ such that $f$ is holomorphic on $\calx \setminus I$. Let $\phi : \Delta \to \calx$ be a complex curve which is not contained in $I$.  Then $\phi$ is holomorphic on $\Delta \setminus \phi^{-1}(I)$. Moreover the subset  $\phi^{-1}(I)$ is discrete.
\smallskip Let $ \lambda_0 \in \phi^{<-1>}(I)$. Consider a ball $B \subset \calx$ centered on $\phi(\lambda_0)$
such that  there exists $F \in \calo( B, l^2)$ with $f|_{B \setminus I} = \pi \circ F$. Then there exists $0<r<1$ such that $ \partial\Delta(\lambda_0 ,r) \cap I = \emptyset$ and $\phi(\Delta(\lambda_0,r)) \subset B$.

\noindent For every $j \in \nn$ we have $F_j \circ \phi (\lambda_0) = 0$. So there exists $k_j \in \nn$ such that $$F_j \circ \phi (\lambda) = (\lambda - \lambda_0)^{k_j} G_j(\lambda)$$ with $G_j \in \calo(B)$ and $G_j$ does not vanish on $\Delta(\lambda_0, r)$. Let $k = \min\{ k_j \ | \ j \in \nn\}$.
Therefore 
$$f \circ \phi(\lambda) = \left[F_0 \circ \phi(\lambda) : F_1 \circ \phi(\lambda) : \dots \right] =  \left[ (\lambda - \lambda_0)^{k_0-k} G_0(\lambda):(\lambda - \lambda_0)^{k_1-k} G_1(\lambda) : \dots \right]$$ with 
$$ (\lambda - \lambda_0)^{2k} \sum_{j=0}^{+\infty} |(\lambda - \lambda_0)^{k_j-k} G_j(\lambda)|^2 = 
\sum_{j=0}^{+\infty} \left|F_j \circ \phi(\lambda) \right|^2 =
\norm{F \circ \phi(\lambda)}^2 < \infty  $$ and then $f \circ \phi$ is holomorphic on $\lambda_0$. 

\newprg[PROJ.ext]{Extension properties of meromorphic maps to $\pp (\l^2)$.}

Now we shall prove the following generalization of the classical Levi, Thullen and 
Hartogs type theorems. We start with Levi-Thullen.

\begin{thm}
\label{levi-thul}
Let $\calx$ be a complex Hilbert manifold and let $f: \calx \setminus A \to \pp(l^2)$ be 
a meromorphic mapping where $A$ is analytic and 

\sli either $\codim A \geqslant 2$ or, more generally,

\slii for every irreducible component $A_{\alpha}$ of $A$ of codimension 
one there exists a point

\quad  $p_{\alpha} \in A_{\alpha}$ such that $f$ meromorphically extends 
to some neighbourhood $V_{\alpha}$ of $p_{\alpha}$. 

\smallskip\noindent Then $f$ meromorphically extends to the whole of $\calx$. 
\end{thm}

\proof \sli We start with codimension $\ge 2$ and denote $A$ as $I$. Our $f$ is defined on 
$\calx\setminus I$ by meromorphic functions $\{f_k\}$ as above. Let $A_0$ 
be as above, \ie $A_0$ is the set of poles of all $f_k$ and their orders along $A_0$ are 
uniformly bounded.  $A_0$ extends to an analytic hypersurface in the whole of $\calx$ by 
Proposition \ref{clos} and $f_k$ meromorphically extend to $\calx$ by Lemma \ref{ext-mero-f}. 
All we need to prove that the orders of poles of these extensions stay uniformly bounded. 
We proceed as in the proof of Proposition \ref{map1}. 
Let $a \in \tilde A_0 \setminus A_0$. There exists a neighbourhood $U$ of $a$ and $g \in 
\calo(U)$  such that $U \cap \tilde A_0  =\{ g=0\}$. Consider the Weierstrass polynomial 
$P$ of $g$ such that $g=P \phi$ where $\phi$ is holomorphic and does not vanish on $U$. 
Since the poles of $f_k$ are uniformly bounded on $U \cap A_0$ we can multiply them by 
$P^{L}$ to obtain a holomorphic map $\psi_k=f_k P^L$ on $U \cap A_0$.  $\psi_k$ extends 
holomorphically to $U$ and then the poles are uniformly bounded on $\tilde{A}_0$.

\smallskip\noindent \slii Represent $f$ by meromorphic functions $\{f_k\}_{k \in \nn^*}$ 
on $(\calx \setminus A)\cup_{\alpha}V_{\alpha}$. By Theorem \ref{thullen} every $f_k$ extends 
meromorphically to the 
whole of $\calx$. Notice that at $p_{\alpha}$ the orders of poles of $f_k$-s are uniformly
bounded by say $N_{\alpha}$. Let $q_{\alpha}\in A_{\alpha}$ be another point. We what to prove 
that the orders of poles of $f_k$-s are bounded at $q_{\alpha}$ with the same $N_{\alpha}$.
Clearly we can assume that both $p_{\alpha}$ and $q_{\alpha}$ are smooth points of $A_{\alpha}$.
Consider a path  $\gamma:[0,1] \to A_{\alpha}^{\reg}$ from $p_{\alpha}$ to $q_{\alpha}$. 
Set
\[
E=\{t\in [0,1] \ | \ \text{ all } f_k \text{ have orders of poles bounded by } N_{\alpha} 
\text{ at } \gamma(t) \}.
\]
$E$ is nonempty and obviously open. Let $t_0= \sup E$ and let $\calv$ be a neighbourhood of 
$\gamma(t_0)$ in $\calx$. Suppose there exists $k$ such that $h^{N_{\alpha}}f_k$ is not holomorphic 
at $\gamma (t_0)$. Take $t_1 <t_0$ such that $t_1\in E$. Now $h^{N_{\alpha}}f_k$ is holomorphic 
in $\calv\setminus A_{\alpha}\cup \{ \text{neighbourhood of } \gamma (t_1)\}$. By Thullen type 
extension for holomorphic functions it extends holomorphically to a neighbourhood of $\gamma (t_0)$.
Contradiction.

\smallskip\qed 

The $n$-dimensional Hartogs figure for $n\ge 2$ is the following open subset of $\cc^n$
\begin{equation}
\eqqno(n-hart)
\calh_1^{n-1}(r)= (\Delta \times B^{n-1}(r)) \cup (A_{1-r,1} \times B^{n-1}).
\end{equation}
Here $r \in ]0,1[$ and $A_{1-r,1}=\Delta \backslash \bar{\Delta}_{1-r} $ is an annulus.
Analogously an infinite dimensional Hartogs figure is the following open subset of $l^2$ 

\begin{equation}
\eqqno(inf-hart)
\calh_1^{\infty}(r)= (\Delta \times B^\infty(r)) \cup (A_{1-r,1} \times B^\infty).
\end{equation}
Denote by $H(r)$ either of these figures and by $B$ the unit ball in $\cc^{n-1}$ or $\l^2$
respectively.
\begin{thm} 
\label{H-ext}
Every $f \in \calm( H(r), \pp(l^2))$ extends meromorphically to $\Delta \times B$.  
\end{thm}
\proof We give the proof in $\l^2$. $f$ is represented by functions $f_k \in \calm(H^\infty_1(r)) 
\cap \calo(H^\infty_1(r) 
\setminus A_0)$. By Corollary \ref{mer-hart} $f_k$ extend meromorphically to $\Delta \times
B^\infty$. Let $A_0$ be the divisor of poles of $f_1$ in $\Delta \times B^\infty$. Every irreducible
component of $A_0$ intersects $\calh_1^{\infty}(r)$. To prove this choose a straight one-dimensional 
closed disk $D\subset B^{\infty}$ in such a way that $\Delta\times D$ intersects the chosen 
component. The question is now reduced to a two-dimensional case.

\smallskip 
Take $a_1 \in A_0^{\reg}$ and choose $a_1  \in \tilde A_0^{\reg}\cap \calh_1^{\infty}(r)$ on the 
same irreducible component as $a_0$. Consider a path $\gamma: [0,1] \to \tilde A_0^{\reg}$ such 
that $\gamma(0)=a_0$ and $\gamma(1)=a_1$. Let $E$ be the set 
\begin{equation}
E =\{t \in [0,1] \ | \text{ all } \ f_k \text{ have bounded denominators in a neighbourhood of } 
\gamma(t)\}.
\end{equation}
Since $f$ is meromorphic on $\calh_1^\infty(r)$ the set $E \neq \emptyset$ and obviously open. 
Denote $t_0 =\sup E$. In a neighbourhood $\calu$ of $\gamma(t_0)$ functions $f_k$ have bounded 
denominators $h_k^N$ where $N$ is the uniform bound. In the same manner as in the proof of Theorem
\ref{levi-thul} one proves that at $\gamma (t_1)$ the orders of poles of $f_k$- are bounded
by $N$.

\smallskip\qed

\smallskip
Recall that a domain $\calx$ in a complex Hilbert manifold with $\calc^2$ boundary
is pseudoconcave at a boundary point 
$a$ if there exists a direction $v \in T^c_a\d \calx$ on which the Levi form of some defining 
function is negative definite. Extension of an analytic object across a pseudoconcave boundary 
point of a domain in a complex Hilbert manifold is equivalent to the Hilbert-Hartogs extension. 
This was explained in \cite{AI} for the case of holomorphic maps. Then the previous theorem 
implies the following

\begin{corol}
\label{concave}
If a domain $\calx$ in a complex Hilbert manifold $\calx$ is pseudoconcave at a boundary 
point $p$ then every meromorphic map $f: \calx \to \mathbb{P}(l^2)$ extends meromorphically 
to a neighbourhood of $p$.
\end{corol}
\proof The problem is local in $\calx$ and therefore we can assume that
$\calx=l^2$. The pseudoconcavity of the point $p \in \partial \calx$ implies that there 
exists $U \subset l^2$ a neighbourhood of $p$ such that $U \cap \calx =\{z \in U \ | \ u(z)<0 \}$ 
with $\nabla u \neq 0$. Take $v\in T_p \partial \calx \subset l^2$ such that  the Levi form 
$\mathcal{L}_{u,p}(v)<0$. After a complex linear change of coordinates we can suppose that 
$p=0$ and $v=e_1,\nabla u_p = e_{2}$.
We define the following figure $ \phi : \calh_1^\infty(r) \to l^2$ by
\begin{equation}
z \mapsto \eta z_1 e_1 + (\eta \epsilon z_{2} -r')e_{2}
+\delta \sum_{s=3}^\infty{ z_se_s}.
\end{equation}

We choose $\epsilon$, $r'$, $\delta$ and $\eta$ such that $\phi(\calh_1^\infty(r)) 
\subset U \cap \calx$ and $\eta \epsilon >r'$ to insure that the image of the polydisk 
$\phi (\Delta \times B^\infty)$ contains the origin. 

Fix $h$ and $j$ and consider now the map $f \circ \phi : \calh_1^\infty(r) \to \pp(l^2)$. By 
Theorem \ref{H-ext} $f\circ \phi$ extends meromorphically to $\tilde{f}: \Delta \times 
B^\infty \to \pp(l^2)$. Then $\tilde{f} \circ \phi^{-1}$ gives the desired extension to 
a neighbourhood of $p$. 

\smallskip\qed

\newsect[LOOP]{Loop space of the Riemann sphere}

\newprg[LOOP.c-str]{Complex structure on the loop space of a complex manifold}

Let $\ss^1$ be the unit circle and consider the set $W^{k,2}(\ss^1, X)$ 
of Sobolev $W^{k,2}$-mappings of $\ss^1$ to a real/complex manifold $X$, 
$\dim_{\rr} X = N\ge 1$ and $N=2n$ with $n=\dim_{\cc} X$ in the complex case.
Notice that since for a smooth function $f$ and an $W^{k,2}$-function $u$ the 
composition $f\circ u$ is again in $W^{k,2}$, see Theorem 1 in \cite{A}, the 
notion of a Sobolev $W^{k,2}$-map from $\ss^1$ to $X$ doesn't depend on the 
coordinate change, \ie is well defined. For the same reason we can  consider 
the pull-back $g^*TX$ of the tangent bundle of $X$ under a $W^{k,2}$-map 
$g:\ss^1\to X$ and get a $W^{k,2}$-bundle.  By that we mean that the 
transition functions of the bundle are in $W^{k,2}$.  Recall furthermore,
that if $g,h\in W^{k,2} (\rr )$ then $gh\in W^{k,2}(\rr)$, and this operation 
is continuous in Sobolev norm, see for example Lemma 2.2 in \cite{A}. This 
enables us to define correctly a $W^{k,2}$-section of our bundle over $\ss^1$.
Moreover, by Sobolev Embedding Theorem we have the following (compact) inclusion 
\begin{equation}
\eqqno(sob-imb)
W^{1,2}(\rr, \rr^N)\subset \calc^{\frac{1}{2}}(\rr,\rr^N).
\end{equation}
In particular, mappings from $W^{1,2}(\ss^1,X)$ are H\"older $\frac{1}{2}$-continuous. 
From now on we shall always assume that $k\ge 1$.
  
\smallskip In a moment we shall see that $W^{k,2}(\ss^1,X)$ is a complex Hilbert 
manifold provided $X$ is complex. But first let us start with its topology.
The Sobolev topology on $W^{k,2}(\ss^1,X)$ can be defined as follows. Take $g\in W^{k,2}
(\ss^1,X)$ and cover $g(\ss^1)$ by finitely many coordinate charts $(X_j,\alpha_j)$, 
$\alpha_j:X_j\to \rr^N$ or $\cc^n$. Similarly cover $\ss^1$ by charts $(S_j,\beta_j)$ 
making sure that $g(\bar S_j)\subset X_j$ for all $j$. In Hilbert spaces $W^{k,2}(S_j,
\rr^N)$ take neighbourhoods $\calv_j$ of $\alpha_j\circ g\circ \beta_j^{-1}$. Then 
define a neighbourhood of $g$ in $W^{k,2}(\ss^1,X)$ as 
\begin{equation}
\eqqno(neib-g)
\calv \deff \{f\in W^{k,2}(\ss^1,X): \alpha_j\circ f\circ \beta_j^{-1}\in \calv_j 
\text{ for all } j\}.
\end{equation}

\smallskip Let us turn to the complex structure of $W^{k,2}(\ss^1,X)$. We closely follow 
\cite{L3}, where the forthcoming lemma was proved in $\calc^{\infty}$ category.
Take $g\in W^{k,2}(\ss^1,X)$. We are going to construct a complex coordinate 
neighbourhood of $g$ in $W^{k,2}(\ss^1,X)$. Consider a neighbourhood $U_g$ of the graph
$\Gamma_g$ of $g$ in $\ss^1\times X$. For every $s\in \ss^1$ consider the following neighbourhood
of $g(s)$ in $X$
\begin{equation}
U_g(s)\deff \{x\in X: (s,x)\in U_g\},
\end{equation}

\begin{lem}
\label{ann-lemp1}
There exists a homeomorphism $G$ between a neighbourhood $U_g$ of $\Gamma_g$ in 
$\ss^1\times X$ and a neighbourhood $V_g$ of the zero section in the induced bundle 
$g^*TX$ such that 

\sli $\{G(s,g(s)):s\in \ss^1\}$ is the zero section of $g^*TX$;

\slii $G_s\deff G|_{U_g(s)}$ maps $U_g(s)$ biholomorphically onto a neighbourhood 
$V_g(s) \deff V_g\cap T_{g(s)}X$ 

\quad of $0\in T_{g(s)}X$.
\end{lem}
Let us recall the main points.
Mapping $G$ is constructed roughly as follows. For a chart $(X_j, \alpha_j)$ in $X$ consider 
a chart $(S_j,\beta_j)$ in $\ss^1$ with $g(S_j)\subset X_j$ as before. Consider an appropriate 
neighbourhood $U_j$ of the graph $g|_{S_j}$ in $S_j\times X_j$ and define $G_j:  U_j\to g^*TX_j$ as
follows
\begin{equation}
\eqqno(clpx-str1)
G_j:(s, x) \to \left(s, (d\alpha_j^{-1})_{\alpha_j(g(s))}[\alpha_j(x) - \alpha_j(g(s))]\right).
\end{equation}
In the right hand of \eqqref(clpx-str1) we use the natural identification  of $T\alpha_j(X_j)$
with $\{\alpha_j(X_j)\}\times \cc^n$ and treat $\alpha_j(x) - \alpha_j(g(s))$ as vector in 
$T_{\alpha_j(g(s))}\alpha_j(X) = \{\alpha_j(g(s))\}\times \cc^n$.
Then glue the maps $G_j$ together to obtain the desired $G(s,x) = \sum_j\eta_j(s)G_j(s,x)$
using a $\calc^k$-partition of unity $\{\eta_j\}$ subordinate to the covering $\{S_j\}$. For 
the more detailed proof  we refer to \cite{L3}, Lemma 2.1.

\smallskip 
Let $U_g$ be a neighbourhood of $\Gamma_g$ in $\ss^1\times X$ as above. Set 
\begin{equation}
\calu_g = \{h\in W^{k,2}(\ss^1,X): \Gamma_h\subset U_g\}.
\end{equation}
For $h\in \calu_g$ set $\Psi_g(h) \deff G(\cdot , h(\cdot))\in W^{k,2}(\ss^1,g^*TX)$. Since $G$ is a homeomorphism, 
$\Psi_g$ is a homeomorphism between $\calu_g$ and a neighbourhood $\calv_g$ of the zero section in 
the space $W^{k,2}(\ss^1,g^*TX)$, which consists from sections that are contained in $V_g$.
We take $(\calu_g,\Psi_g)$ as a complex $\l^2$- coordinate chart in $W^{k,2}(\ss^1, X)$. Let 
$h\in W^{k,2}(\ss^1,X)$ be such that $\Gamma_h\subset U_g\cap U_{g'}$, 
\ie $h\in \calu_g\cap\calu_{g'}$. For $s\in \ss^1$ we have 
\[
\Psi_{g'}(h)(s) = \left(\Psi_{g'}\circ \Psi_g^{-1}\right)\Psi_g(h)(s) = 
\left(G'_s\circ G_s^{-1}\right)\left(\Psi_g(h)(s)\right) = 
\]
\begin{equation}
= \left(G'_s\circ G_s^{-1}\right) \left(\Psi_g(h)(s)\right).
\end{equation}

\smallskip Due to the item (\slii of Lemma \ref{ann-lemp1} for every $s$ mapping 
$G'_s\circ G_s^{-1}$ is a biholomorphism between an appropriate open 
subsets of $g^*T_{g(s)}X$ and $g'^*T_{g'(s)}X$. Let us prove that \( \Psi_{g'}\circ \Psi_g^{-1} \) is holomorphic. Since it is a homeomorphism, it is in particular continuous. We shall now prove that it is Gâteaux holomorphic.

Let $ v, w \in W^{k,2}(\ss^1, g^*TX) $ and $ \lambda \in \Delta(0,\epsilon) $ with $ \epsilon $ small enough. For every $ s \in \ss^1 $, the holomorphicity of $ G'_s\circ G_s^{-1} $ gives us

\begin{equation}
\Psi_{g'} \circ \Psi_g^{-1} (v + \lambda w)(s) = G'_s(\Psi_g^{-1}(v + \lambda w)(s)) = 
G'_s( (G_s)^{-1}(v(s) + \lambda w(s))) = 
\end{equation}
\[
= \frac{1}{2\pi i }\int_{\Delta(0,\epsilon)} \frac{1}{\xi - \lambda}G'_s\big((G_s)^{-1}(v(s) + 
\xi w(s))\big) d\xi = 
\]
\[
 =  \left( \frac{1}{2\pi i }\int_{\Delta(0,\eps)} \frac{1}{\xi - \lambda}\Psi_{g'} \circ 
\Psi_g^{-1} (v + \xi w) d\xi \right) (s).
\]
Since this holds for every $ s \in \ss^1 $, we obtain
\begin{equation}
	\Psi_{g'} \circ \Psi_g^{-1} (v+ \lambda w) = \frac{1}{2\pi i }\int_{\Delta(0,\epsilon)} \frac{1}{\xi - \lambda}\Psi_{g'} \circ \Psi_g^{-1} (v+ \xi w) d\xi.
\end{equation}
Therefore $ \Psi_{g'} \circ \Psi_g^{-1} : W^{k,2}(\ss^1 , g^*TX) \to W^{k,2}(\ss^1 , g'^*TX) $ is Gâteaux holomorphic and continuous. It follows that it is holomorphic and therefore is a biholomorphism.

\smallskip From this description of the complex structure in $W^{k,2} (\ss^1,X)$
one can deduce, see \cite{L1} in $\calc^\infty$ category, or Proposition 2.3 and
Lemma 3.2 in \cite{A} for case of Sobolev regularity, the following:
\begin{prop}
\label{lempert}
A holomorphic map from a complex Hilbert manifold $\calx$ to $W^{k,2}(\ss^1, X)$ can 
be represented by a mapping $f:\calx \times \ss^1\to X$ such that:

\smallskip\sli for every $s\in \ss^1$ mapping $f(\cdot,s):\calx\to X$ is holomorphic;

\slii for every $z\in \calx$ one has $f(z,\cdot)\in W^{k,2} (\ss^1,X)$ and the
correspondence 
\[
\calx\ni z\to f(z,\cdot)\in W^{k,2} (\ss^1,X)
\]

\quad is continuous with respect to the Sobolev topology on $W^{k,2}(\ss^1,X)$ and the  

\quad standard topology on $\calx$. Mapping $f:\calx \to W^{k,2}(\ss^1,X)$ is actually 
$z\to f(z, \cdot)$.
\end{prop}

\begin{rema} \rm
\label{ann-lemp2}
{\bf a)} Notice that if $X=\cc^n$ then $G(s,x) = (s,x -g(s))$ because $TX\equiv X\times \cc^n$ in this 
case. The same simple form has $G$ if the image of $g$ is contained in one chart $X_j\sim \tilde X_j
\subset \cc^n$. Therefore $\Psi_g$ has the form $\Psi_g(h)(s) = (s, h(s) - g(s))$ in these cases.

\smallskip\noindent
{\bf b)} And one more remark, let $f$ be a holomorphic function on $X_j$ and let $\calx_j$ be the
open set of such $h\in W^{k,2}(\ss^1,X)$ that $h(s_0)\in X_j$. Then $F(h)\deff f(h(s_0))$ is a 
holomorphic function on $\calx_j$. Indeed, we can assume that $X_j$ is a coordinate neighbourhood
of some $g(\ss^1)$ and consider loops $h$ close to $g$. Then using the coordinate chart 
$(\calx_j , G)$ as above we see that $F(h) = f\left((G^{s_0})^{-1}(h(s_0)-g(s_0))\right)$ in the 
notation of Lemma \ref{ann-lemp1}. This proves the holomorphicity of $F$.

\smallskip\noindent{\bf c)} If $f$ is meromorphic write (locally) $f=p/q$ and set $F(g) = 
f(h(s_0)) = p(h(s_0))/q(h(s_0))$ in the case when $h(s_0)\not\in I_f$. This will be a meromorphic 
function on $\calx_j\setminus \{h\in\calx_j: h(s_0)\in I_f\}$. The analytic set $\{h\in\calx_j: 
h(s_0)\in I_f\}$ is of codimension $\ge 2$ in $\calx_j$ and therefore $F$ extends to a meromorphic
function on $\calx_j$. If $f$ was meromorphic on the whole of $X$ we obtain a meromorphic function 
$F$ on the whole of $W^{k,2}(\ss^1, X)$, in fact a family of meromorphic functions parametrized by
$s_0\in \ss^1$.

\smallskip\noindent
{\bf d)} Finally let us remark that $L\pp^1$ is homogeneous since the loop group $LG$, where 
$G=PGL(2,\cc)$ obviously acts holomorphically and transitively on $L\pp^1$.
\end{rema}

\begin{prop}
\label{mer-sep}
Let $X$ be a complex manifold such that $\calm (X)$ separates points on $X$. Then 
$\calm (W^{k,2}(\ss^1,X))$ separates points on $W^{k,2}(\ss^1,X)$.
\end{prop}
For if $g,h\in W^{k,2}(\ss^1,X)$ are distinct then $g(s_0)\not=h(s_0)$ for some $s_0\in \ss^1$.
By assumption there exists a meromorphic function on $X$ which separates $g(s_0)$ and $h(s_0)$. 
The corresponding $F$ constructed as in the remark above will separate $g$ and $h$.

\smallskip We finish this subsection with statement mentioned in the Introduction.

\begin{prop}
\label{mer-inj}
Let $X$ be a finite dimensional compact complex manifold such that $\calm (X)$ separates 
points on $X$. Then $X$ is Moishezon and therefore admits a meromorphic injection to the
complex projective space $\pp^N$ for some $N$.
\end{prop}
\proof Let $I_f$ be the indeterminacy set of $f$. The proof will be done in two steps. 

\smallskip\noindent{\slsf Step 1. {\sl $X$ is Moishezon.}} Suppose not, 
\ie $t\deff \trdeg \calm (X) < n = \dim X$. Take some transcendence base $g_1,...,g_t$
of $\calm (X)$. Then for a generic point $c=(c_1,...,c_t)
\in \cc^t$ there exists a connected component $F_c$ of the level set 
$\{g_1=c_1,...,g_t=c_t\}$ of dimension $n - t$ such that any meromorphic 
function cannot separate points in $F_c$. Indeed, if we suppose that there exists 
meromorphic function $f$, which is regular near $x\not= y\in F_c$, and 
such that $f(x)\not=f(y)$, then $f$ satisfies the equation 
\begin{equation}
\eqqno(eq-f)
f^N+h_1(g_1,...,g_t)f^{N-1}+...+h_N(g_1,...,g_t)=0, 
\end{equation}
where $h_i$ are rational functions. It may happen that $x$ or $y$ belong to the 
poles of $h_i(g_j)$. But then we can perturb $c$ in such a way that the 
deformed connected component $F_c$ is not contained entirely in the 
union of these poles. From \eqqref(eq-f) we see  that $f=\const $ 
along $F_c$. Contradiction. 

\smallskip\noindent{\slsf Step 2. {\sl Construction of a meromorphic injection.}}
Since $X$ is Moishezon, there exists a bimeromorphic map $f:X\merto Y$, where 
$Y$ is projective. Apply the Hironaka resolution of singularities theorem to 
the inverse map $f^{-1}:Y\merto X$ and get a holomorphic map $\tilde f^{-1} : 
\tilde Y \to X$ which is bimeromorphic. Here $\tilde Y$ is obtained by 
a finite sequence of blowings up of $Y$ with smooth centers, see \cite{BM}. 
This $\tilde Y$ is projective. Take the inverse map to $\tilde f^{-1}$, 
denote it as $\tilde f:X\to \tilde Y$. This map doesn't contract 
divisors in $X$ (because its inverse is regular) and is bimeromorphic.

\smallskip Suppose its differential vanishes at some point $x\in X\setminus
I_{\tilde f}$. If $x$ is taken generic then $\tilde f$ is a ramified 
cover in a neighbourhood of $x$. But this contradicts to the fact that 
$\tilde f$ is bimeromorphic. Therefore $\tilde f$ is an embedding on 
$X\setminus I_{\tilde f}$.

\smallskip\qed

\newprg[LOOP.mero]{Meromorphic maps to the loop space}
From now on we assume that $X=\pp^1$, $k=1$ and we are going to explain that the Definition 
\ref{mer-map2} from the Introduction gives the natural notion of a meromorphic mapping 
in the case of $L\pp^1=W^{1,2}(\ss^1,\pp^1)$ as a target manifold. Let a meromorphic 
(in any sense) mapping $f: \calx \to L\pp^1$ be given. Then the least thing we know is 
that there should exist 
an analytic of codimension $\geqslant 2$ set $I$ in $\calx$  such that the restriction 
of $f$ to $\calx \setminus I$ is holomorphic. Denote this restriction still as $f:\calx
\setminus I \to L\pp^1$ as well as its representation $f:(\calx\setminus I)\times \ss^1
\to \pp^1$. Notice that the latter satisfies items (\sli and (\slii of Proposition 
\ref{lempert}. 

\begin{thm}
\label{loop-mer1}
Let a holomorphic mapping $f:\calx\setminus I\to L\pp^1$ be given, where $I$ is analytic 
in $\calx$ of codimension $\ge 2$. Then:

\sli for every $s\in \ss^1$ mapping $f(\cdot , s)$ is a meromorphic function which extends 
mero-

\quad  morphically on the whole of $\calx$ and its indeterminacy set $I_{f(\cdot , s)}$
is contained in $I$;

\slii the family of the indeterminacy sets $\{I_{f(\cdot , s)}:s\in \ss^1\}$ is locally
finite in $\calx$ and 

\quad $I_f = \bigcup_{s \in \ss^1}I_{f(\cdot, s)}$ is an analytic set in $\calx$ of pure 
codimension two;

\sliii mapping $f$ holomorphically extends to $\calx\setminus I_f$.
\end{thm}
\proof The proof will be achieved in several steps.

\smallskip\noindent{\slsf Step 1.} {\it Meromorphicity of $f(\cdot ,s)$.} Notice that 
for every $s \in \ss^1$ mapping $f(\cdot, s):\calx\setminus I\to \pp^1$ is in fact a 
meromorphic function without indeterminacy points and it extends to a meromorphic 
function on the whole of $\calx$ by Lemma \ref{ext-mero-f}. Denote by $I_{f(\cdot , s)}$ 
the indeterminacy set of the extended $f(\cdot ,s)$. We have that $I_{f(\cdot ,s)}\subset I$ 
because $f$ is holomorphic on $\calx\setminus I$. Since $I_{f(\cdot ,s)}$ is the 
indeterminacy set of a meromorphic function and, therefore is locally an intersection of 
two relatively prime divisors it is of pure codimension two. 
Take some point $z_0\in \bigcup_{s\in \ss^1}I_{f(\cdot ,s)}\subset I^2_{z_0}$, where 
$I^2_{z_0}$ is the (necessarily finite) union of codimension two components of $I$ at 
$z_0$ as in \eqqref(decomp1). Due to the description of the structure of Hilbert analytic 
sets of finite definition in \cite{Ra} one can find a neighbourhood $\calu = B^{\infty}\times 
\Delta^2$ of $z_0$ such that $I^2_{z_0}\cap \calu$ is a finite cover of $B^{\infty}$ in these 
coordinates. Therefore $\left(\bigcup_{s\in \ss^1}I_{f(\cdot ,s)}\right)\cap \calu$ can 
have only finitely many components. We obtain as a result that the union $\bigcup_{s \in 
\ss^1} I_{f(\cdot, s)}$ is an analytic of pure codimension two subset of $\calx$.

\smallskip Items (\sli and (\slii are proved and we turn to item (\sliii. Take some $z_0 
\in \calx\setminus \bigcup_{s \in \ss^1} I_{f(\cdot, s)}$ and let $\calu$ be a neighbourhood 
of $z_0$ not intersecting $\bigcup_{s \in \ss^1} I_{f(\cdot, s)}$. Then $f$ extends to a 
mapping $f:\calu\times \ss^1\to \pp^1$. We need to prove its holomorphicity. Condition (\sli 
of Proposition \ref{lempert} for our $f$ is obviously satisfied. Therefore all we need 
to prove is the condition  (\slii , \ie continuity of the mapping $z\to f(z, \cdot)$ in 
Sobolev topology. 

\smallskip\noindent{\slsf Step 2.} {\it  Continuity in classical topology.} We shall first prove
that $f:\calu\times \ss^1\to\pp^1$ is continuous in the classical sense. We need to prove this at the
points $(z_0,s_0)\in \calu\times \ss^1$ with $z_0\in I\setminus I_f$ only. Fix such point $z_0$ 
and take an open neighbourhood $\calu$ of $z_0$ biholomorphic to $B^{\infty}\times\Delta^2$, where 
bidisk $\{0\}\times \Delta^2$ intersects $I$ at $z_0$ transversely, \ie $z_0$ is an isolated point
of $(\{0\}\times \Delta^2)\cap I$. Natural coordinates in $B^{\infty}\times\Delta^2$ centered at 
$z_0 = (0,0)$ denote as $z=(w,t)$. Note  that $f$ is holomorphic on $B^{\infty}\times\d\Delta^2$, 
provided $B^{\infty}$ was taken sufficiently small, and that $f(w,\cdot , s_0)|_{\Delta^2}$ is 
holomorphic over the whole of $\{w\}\times\Delta^2$ for every fixed $w\in B^{\infty}$. Due to 
Lemma 3.1 from \cite{AZ} to prove the continuity of $f$ at $(z_0,s_0)$ it is sufficient to prove 
that for any sequences $s_n\to s_0$ in $\ss^1$, $w_n\to 0$ in $B^{\infty}$ and any neighbourhood 
$\calv$ of the graph $\Gamma_0$ of $f(0, \cdot , s_0)|_{\Delta^2}$ in $\calu\times \pp^1$ the 
graphs $\Gamma_n$ of $f(w_n, \cdot , s_n)|_{\Delta^2}$ are contained in $\calv$ for $n$ big 
enough. By the main result of \cite{AZ}, \ie Theorem 2.1 there, we can assume that $\calv$ is 
$1$-complete. That is it admits a strictly plurisubharmonic exhausting function $\rho$. Now it 
is easy to conclude from Sobolev embedding \eqqref(sob-imb) that for $w_n$ and $s_n$ close enough 
to $0$ and $s_0$ respectively $f(w_n, \cdot ,s_n)|_{\d\Delta^2}$ is close to 
$f(0, \cdot ,s_0)|_{\d\Delta^2}$. And than by maximum principle for $\rho$ to deduce that 
$\Gamma_n$ stays in $\calv$. The step is proved.

\smallskip\noindent{\slsf Step 3.} {\it Continuity in Sobolev topology.} Now we need to prove
that a continuous in classical sense mapping $f:\calu\times \ss^1\to \pp^1$ is continuous also 
in Sobolev topology. Using compactness of $\ss^1$ cover it by arcs $S_j(\eps) =]s_j-\eps , s_j+
\eps [ , j=1,...,N$ centered at $s_j$ and shrink a neighbourhood $\calu$ of $z_0$ so that 
$f(\calu\times S_j)$ is contained in a coordinate chart $(\calv_j,\alpha_j)$ of $\pp^1$. Make 
sure that $f\left(\overline{\calu\times S_j(\eps_1)}\right) \subset \calv_j$ still holds true 
for some $\eps_1>\eps$. Consider the mapping $f_j\deff \rho_j\alpha_j\circ f|_{\calu\times 
S_j(\eps_1)}$, where $\rho_j$ is a cut-off function supported in $S_j(\eps_1)$ and equal to 
$\adyn$ on $S_j(\eps)$. Take $\calu$ this time in the form $\Delta\times B^{\infty}$ and find 
a Hartogs figure $H^{\infty}_1(r)=\Delta\times B^{\infty}(r)\cup A_{1-r,1}\times B^{\infty}
\subset \Delta\times B^{\infty}$ not intersecting $I$. Let $(z,t)$ be the natural coordinates 
in $\Delta\times B^{\infty}$.

\smallskip Note that: 

\smallskip a) for every $(z,t)\in H^{\infty}_1(r)$ mapping $f_j$ is in $W^{1,2}(S_j(\eps_1), \cc)$
and has compact support.

\smallskip b) Mapping $(z,t) \to f_j(z,t,\cdot)\in W^{1,2}(S_j(\eps_1),\cc)$ is continuous in Sobolev 
topology.

\smallskip c) Moreover, for every $s\in S_j(\eps_1)$ mapping $f_j(\cdot , s)$ is holomorphic on 
$\hat H^{\infty}_1(r) = \Delta\times B^{\infty}$.

\smallskip The latter is because the multiplication by a smooth cutoff function preserves the 
Sobolev class and continuity in Sobolev topology and since it depends only on the space variable
$s$ it doesn't spoils the holomorphicity in $(z,t)$. Therefore $f_j$ is a holomorphic mapping 
from $H^{\infty}_1(r)$ to a complex Hilbert space. Consequently it extends holomorphically to 
$\Delta\times B^{\infty}$. In particular we obtain that $f|_{\calu\times S_j(\eps)}$ is 
continuous in Sobolev topology. By the definition of the Sobolev topology, see \eqqref(neib-g) 
we obtain the desired continuity of $f$. Lemma is proved.

\smallskip\qed 

\begin{rema} \rm
\label{def-mer-l}
As the result we reach the understanding that the Definition \ref{mer-map2} from the Introduction 
is the natural definition of meromorphicity for mappings with values in $L\pp^1$. Indeed, what
we proved in Proposition \ref{loop-mer1} is that a holomorphic map $f:\calx\setminus I\to L\pp^1$,
where $I$ is Hilbert analytic of codimension $\ge 2$, extends to a meromorphic map from 
$\calx$ to $L\pp^1$ in the sense of Definition \ref{mer-map2}.
\end{rema}

\newprg[LOOP.imb]{Non-existence of an embedding of the loop space to the projective space}
We start this subsection with the following example.
\begin{exmp} \rm
\label{loop-exmp}
Let $\gamma (s)= (\gamma_1(s),\gamma_2(s))$ be a loop in $\cc^2$ of class $W^{1,2}$. Denote 
by $\Gamma$ the image of $\gamma$. Consider the mapping $f: \cc^2\setminus \Gamma \to 
W^{1,2}(\ss^1,\pp^1)$ defined as follows
\begin{equation}
\eqqno(exmpLoop)
 f(z_1,z_2,s) = [z_1-\gamma_1(s):z_2-\gamma_2(s)].
\end{equation}
For every fixed $s\in \ss^1$ mapping $f(\cdot ,s):\cc^2\setminus \Gamma \to \pp^1$ is  
holomorphic. Moreover, for every $z\in \cc^2\setminus \Gamma$ we have that $f(z,\cdot) \in 
W^{1,2} (\ss^1,\pp^1)$ and it continuously depends on $z$. This means that $f$ is a holomorphic
mapping from $\cc^2\setminus \Gamma$ to $W^{1,2}(\ss^1,\pp^1)$. Furthermore mapping 
$f(\cdot ,s)$ is in fact holomorphic on $\cc^2\setminus \{\gamma (s)\}$, having $\gamma(s)$ as its
only indeterminacy point. Therefore  $\bigcup_{s \in \ss^1} I_{f(\cdot, s)}=\Gamma$, \ie 
the family $\{I_{f(\cdot ,s)}\}$ is not locally finite and henceforth is not an analytic set in 
$\cc^2$. We see that $\Gamma$ is an {\sl essential singularity} of $f$ in the sense that $f$ 
does not extend meromorphically to a neighbourhood of any point of $\Gamma$. Indeed, if $f$ 
extends meromorphically to a neighbourhood $V$ of some $\gamma(s_0) \in \Gamma$ it should be 
holomorphic on $V$ minus a discrete set. But $f$ is obviously not holomorphic at any point 
of $\Gamma \cap V$.  
\end{exmp}

\begin{rema} \rm 
\label{non-thul}
This example shows that the Hartogs type extension theorem is not valid for 
meromorphic values in $L\pp^1$. Taking $\Gamma \subset \cc \subset\cc^2$ we 
obtain a counter example also for the Thullen type extension.
\end{rema}

\medskip\noindent{\slsf Proof of Theorem \ref{non-imb}.}
Suppose to the contrary that $g:L\pp^1 \merto \pp(l^2)$ is such a map and let $\gamma$ be such
as in Example \ref{loop-exmp}. Take $\gamma =(\gamma_1,\gamma_2)$  such that it is of class 
$W^{2,2}\subset \calc^{1,\frac{1}{2}}$, both components are an immersed curves and in addition
that its image $\Gamma$ is polynomially convex. In particular 
it is {\slsf not contained} in any complex curve in $\cc^2$. Consider the holomorphic map 
$f: \cc^2\setminus \Gamma \to L\pp^1$ defined in \eqqref(exmpLoop). Using the homogeneity of 
$L\pp^1$ compose $f$ with an automorphism of $L\pp^1$ in order to map some point $p\in 
\cc^2\setminus \Gamma$ by a resulting map to a  point $q\in L\pp^1$ of $g$ such 
that $\codim\ker dg_q\ge 2$. Denote this resulting map still as $f$.
\begin{lem}
\label{inj}
The differential of mapping $f:\cc^2\setminus \Gamma\to L\pp^1$ is everywhere injective.
\end{lem}
\proof We need to prove this for the original $f$ as in \eqqref(exmpLoop) since 
a composition with an automorphism doesn't affect the injectivity. Take $z^0=(z_1^0,z_2^0)\in 
\cc\setminus \Gamma$ such that $z_1^0-\gamma_1(s)$ never vanishes. Then 
\[
f(z,s) = [z_1-\gamma_1(s):z_2-\gamma_2(s)]
\]
takes its values in $U_0 = \{[w]=[w_0:w_1]\in \pp^1: w_0\not=0\}$ for $z$ in some neighbourhood 
$U\ni z^0$. I.e., 
\[
f(z,s)=\frac{z_2-\gamma_2(s)}{z_1-\gamma_1(s)}
\]
represents a holomorphic map $f:U\to W^{1,2}(\ss^1,\cc)$. Now for $z\in U$ and $\vect \in \cc^2$
one has 
\begin{equation}
\eqqno(diff-f1)
df_z[\vect] = \frac{\d f(z,s)}{\d z_1}\vect_1+\frac{\d f(z,s)}{\d z_2}\vect_2 = \left(s, 
\frac{z_2-\gamma_2(s)}{[z_1-\gamma_1(s)]^2}\vect_1 + \frac{\vect_2}{z_1-\gamma_1(s)}\right).
\end{equation}
Here the right hand side of \eqqref(diff-f1) should be understood as a section of the pulled-back 
tangent bundle $f^*(z,\cdot)T\pp^1\equiv \ss^1\times \cc$. To prove the injectivity of $df_z$ 
assume that $df_z[\vect] = 0$ for some $\vect\not= 0$. Then 
\[
\frac{z_2-\gamma_2(s)}{[z_1-\gamma_1(s)]^2}\vect_1 + \frac{\vect_2}{z_1-\gamma_1(s)} = 0
\]
for all $s$. And this means that 
\[
\frac{z_2-\gamma_2(s)}{z_1-\gamma_1(s)}\vect_1 + \vect_2 = 0.
\]
Since $\vect\not= 0$ this implies that $\frac{z_2-\gamma_2(s)}{z_1-\gamma_1(s)} = \const$ and 
therefore $\gamma_2(s) - \c_1 \gamma_1(s) = \c_2$ for some constants $\c_1, \c_2$ and all $s$. 
I.e., $\Gamma$ is contained in the complex line $C = \{z_2\c_1-z_1 = \c_2\}$ and henceforth 
cannot be polynomially convex. Contradiction. Analogously one proves that $df_z$ is injective 
in an neighbourhood of any point $z^0=(z_1^0,z_2^0)$ such that $z_2^0-\gamma_2(s)$ never vanishes.

\smallskip We proved that $df_z$ can degenerate for $z$ in a complex curve $C\subset 
\left(\gamma_1(\ss^1)\times \gamma_2(\ss^1)\right)\setminus \Gamma $ at most. But 
$\gamma_1(\ss^1)\times \gamma_2(\ss^1)$ is totally real and therefore $C=\emptyset$.

\smallskip\qed 

\begin{lem}
\label{transv}
After an appropriate perturbation $\tilde \gamma$ of $\gamma$ the map $\tilde f$
thus obtained will satisfy $df_p[\cc^2]\pitchfork \ker dg_q$.
\end{lem}
\proof Injectivity of $df_p$ together with \eqqref(diff-f1) show that 
tangent vectors 
\begin{equation}
\eqqno(transv)
\begin{cases}
 w_1 = p_2-\gamma_2(s) = df_p[e_1]\cr 
 w_2 = p_1-\gamma_1(s) = df_p[e_2]
\end{cases}
\end{equation}
are linearly independent and generate $df_p[\cc^2]$, here $p=(p_1,p_2)$. If one, or both of 
$w_1,w_2$  belong to $\ker dg_q$ then we can perturb them in order become transverse 
to $\ker dg_q$. And then appropriately perturbe $\gamma_1,\gamma_2$ for \eqqref(transv)
to stay valid. All what is left is to remark that small perturbations of $\gamma_k$
will stay to be immersed and a small perturbation of $\Gamma =\gamma (\ss^1)$ will stay
to be polynomially convex.

\smallskip\qed

Now the composition $h\deff g\circ f$ is well defined and 
by Corollary \ref{hart-corol} the map $h: \cc^2 \setminus \Gamma \to \pp(l^2)$ extends
meromorphically onto $\cc^2$. Denote this extension still as $h$. Take $a \in \Gamma 
\setminus I_h$ and consider a neighbourhood $V$ of $a$ which do not contain points from 
$I_h$, \ie the restriction $h|_V: V \to \pp(l^2)$ is holomorphic. Notice that its differential 
is injective at all points of $V$.

\smallskip Furthermore, $dh_z[\cc^2]$ is closed in $T_{h(z)}\pp (\l^2)$ for all $z\in V$ since 
it is a two-dimensional subspace of $T_{h(z)}\pp (\l^2)\equiv \l^2$. We find an affine chart $\calu$ 
containing $h(V)$. For this one might need to shrink $V$ again. Now we compose $h$ with the 
orthogonal projection $\pi$ of $\l^2$ to $dh_a[\cc^2]\equiv \cc^2$.  This composition is 
biholomorphic in a neighbourhood (again call it $V$) of $a$. Therefore $h(V)$ is a graph 
over $\pi (h(V))$ in $\cc^2\times \l^2$. This proves that $h|_V:V\to h(V)$ is a biholomorphism
to a local submanifold $h(V)$ of $g(L\pp^1)$. $g^{-1}$ is well defined on $h(V)$ and therefore 
$g^{-1}\circ h$ is holomorphic on $V$. But it coincides with $f$ on $V\setminus \Gamma $. We got
a holomorphic extension of $f$ from $V\setminus \Gamma$ to $V$. Composing this extension with 
the inverse to the automorphism of $L\pp^1$ taken at the beginning of the proof we get an 
extension of the original $f$. This contradiction proves the Theorem.

\smallskip\qed

\ifx\undefined\bysame
\newcommand{\bysame}{\leavevmode\hbox to3em{\hrulefill}\,}
\fi

\def\entry#1#2#3#4\par{\bibitem[#1]{#1}
{\textsc{#2 }}{\sl{#3} }#4\par\vskip2pt}

\end{document}